\documentclass[12pt,a4paper]{amsart}
\usepackage[T1]{fontenc}
\usepackage[utf8]{inputenc}
\usepackage{amssymb, amsthm, latexsym}
\usepackage[margin=3cm]{geometry}

\theoremstyle{plain}
\newtheorem{thm}{Theorem}
\newtheorem{lemma}[thm]{Lemma}
\newtheorem{corollary}[thm]{Corollary}
\newtheorem{proposition}[thm]{Proposition}
\newtheorem{fact}[thm]{Fact}

\theoremstyle{definition}
\newtheorem{definition}[thm]{Definition}
\newtheorem{remark}[thm]{Remark}
\newtheorem{example}[thm]{Example}

\numberwithin{thm}{section}
\numberwithin{equation}{section}

\renewcommand{\a}{\alpha}
\newcommand{\e}{\varepsilon}
\renewcommand{\d}{\delta}
\renewcommand{\k}{\kappa}
\renewcommand{\l}{\lambda}
\newcommand{\n}{\eta}
\renewcommand{\o}{\omega}
\newcommand{\distp}{\textup{d}^p}
\newcommand{\raj}{\restriction}
\newcommand{\A}{\mathcal{A}}
\newcommand{\B}{\mathcal{B}}

\newcommand{\D}{\mathcal{D}}
\newcommand{\F}{\mathbf{F}}
\newcommand{\K}{\mathbf{K}}
\newcommand{\M}{\mathbf{M}}
\newcommand{\N}{\mathbb{N}}
\newcommand{\Q}{\mathbb{Q}}
\newcommand{\Z}{\mathbb{Z}}
\newcommand{\C}{\mathbb{C}}
\newcommand{\R}{\mathbb{R}}
\renewcommand{\empty}{\emptyset}
\newcommand{\abs}[1]{\lvert#1\rvert}
\newcommand{\norm}[1]{\lVert#1\rVert}

\newcommand{\tg}{t^g}
\newcommand{\BF}{\mathbb{F}}
\newcommand{\Mon}{\M} 
\newcommand{\de}{\operatorname{d}}

\def\Ind#1#2{#1\setbox0=\hbox{$#1x$}\kern\wd0\hbox to 0pt{\hss$#1\mid$\hss}
\lower.9\ht0\hbox to 0pt{\hss$#1\smile$\hss}\kern\wd0}
\def\ind{\mathop{\mathpalette\Ind{}}}
\def\notind#1#2{#1\setbox0=\hbox{$#1x$}\kern\wd0
\hbox to 0pt{\mathchardef\nn=12854\hss$#1\nn$\kern1.4\wd0\hss}
\hbox to 0pt{\hss$#1\mid$\hss}\lower.9\ht0 \hbox to 0pt{\hss$#1\smile$\hss}\kern\wd0}
\def\nind{\mathop{\mathpalette\notind{}}}
\def\da{\ind}
\def\nda{\nind}

\begin{document}

\title[Measuring dependence in MAECs with perturbations]{Measuring dependence in metric abstract elementary classes with perturbations}
\date{\today}
\author{{\AA}sa Hirvonen}
\address{Department of Mathematics and Statistics\\
University of Helsinki\\
P.O.Box 68\\
00014 University of Helsinki\\
Finland}
\email{asa.hirvonen@helsinki.fi}
\author{Tapani Hyttinen}
\address{Department of Mathematics and Statistics\\
University of Helsinki\\
P.O.Box 68\\
00014 University of Helsinki\\
Finland}
\email{tapani.hyttinen@helsinki.fi}
\thanks{Both authors were partially supported by Academy of Finland grant 251557. The first author was also supported by the Ruth and Nils-Erik Stenb\"ack foundation.}
\begin{abstract}
We define and study a metric independence notion in a
homogeneous metric abstract elementary class with perturbations
that is $d^p$-superstable (superstable wrt.\ the perturbation topology),
weakly simple and has complete type spaces
and we give a new example of such a class based on B. Zilber's
approximations of  Weyl algebras.
We introduce a way to measure the dependence of a tuple $a$ from
a set $B$ over another set $A$. We prove basic properties of the notion,
e.g.\ that $a$ is independent of $B$ over $A$
in the usual sense of homogeneous model theory if and only if the measure
of dependence is $<\e$
for all $\e >0$.
As an example of our measure of dependence we show a connection between the
measure and entropy in  models from quantum mechanics in which
the spectrum of the observable is discrete.
As an application, we show
that weak simplicity implies a very strong form of simplicity
and study the question of when the dependence inside a set of
all realisations of some type can be seen to arise from
a pregeometry in cases when the type is not regular.
In the end of the paper, we demonstrate our notions and results
in
one more example: a
class built from the $p$-adic integers.
\end{abstract}

\maketitle

\section{Introduction}

Studying metric structures from a stability theoretic point of view
offers two basic new features. The first is to be expected: the
natural notion of size has to switch from cardinality to density. The
second was noted by Iovino \cite{iovstab} and developed further by Ben
Yaacov \cite{benpert}: studying stable classes of metric structures,
the cardinalities at which stability occurs depends on the metric (or
more generally, topology)
chosen for the type space. Ben Yaacov's approach built on the idea of
allowing \emph{perturbations} to the structures, launching the term
\emph{$\omega$-stable up to perturbation}. Examples of structures
which are $\omega$-stable up to perturbation (but not as nonperturbed)
were studied in \cite{buz} and \cite{byberpert}.

To be able to use the improvement in stablility we introduced the
framework of \emph{metric abstract elementary classes with
  perturbations}  in \cite{a2}. The idea is to work in a syntax-free
framework and add the perturbation mappings as graded classes of
generalised isomorphisms. As the main motivation for the approach is
to enable a natural use of perturbation mappings, but not to work in
as general a setting as conceivable (as opposed to the general framework of
abstract elementary classes), we have worked in a homogeneous
context. So, were it not for the perturbations, the setting would roughly correspond to a metric variant of Shelah's finite diagrams, and the reader may think of the types as the syntactic types of a homogeneous monster. However, the topology on the type space need not be given by any language. A first development of basic independence
and isolation notions under $\omega$-$\distp$-stability ($\omega$-stability up to perturbations) in this context was presented in \cite{a2}.

This paper continues the study of homogeneous metric structures with
perturbations. Stability requirements are now reduced to
$\distp$-superstability (together with some other requirements, specified below), but the
main novelty is not just this loosening of stability requirements
considered, but the way the independence notion developed
offers a measure of the dependence.

We define $\distp$-superstability as `stable with respect to the $\distp$-topology from some $\lambda$ onward' and show that it is equivalent to other natural characterisations. As far as generalising stability assumptions is concerned, assuming $\distp$-superstability is as far as is reasonable to go, as perturbations do not bring stability to unstable classes, but only brings down the stability spectrum. So assuming only $\distp$-stability, one would be tied to the effort of dealing with perturbations, without gaining access to better methods than what could already be used by studying the class, say, as a discrete homogeneous AEC\@. In the context of AEC's,
the main point of looking at a class as a metric class is that
this may bring the class down in the stability hierarchy allowing for better theory.

We also assume completeness of type spaces. This is a weakening of the compactness property one has in continuous first order logic, so the assumption is satisfied if the class consists of all the models of some
complete theory in continuous logic and the perturbations are either trivial or can be captured by approximating formulas (a case, which the reader may use as a first mental picture). But there are more general classes satisfying the assumptions, e.g., the class of real valued atomic Nakano spaces (treated in \cite{a2}).

We aim to create a theory for
measuring dependencies. Already in \cite{dfin}
something like this was done but
here our approach is a bit different and we aim to go further.
The main difficulty in developing the theory is that independence
is closely linked with
properties of
Lascar (strong) types.
So to develop the theory one should be
able to measure the distances between Lascar types. It is not clear how this could be done in general, however, if the class is simple there is a way. Here by simple we
mean that all Galois-types over any set have a free extension
to any other set. Since our independence notion will be defined with built-in
free extensions, this is the same as saying that for all $a$ and $A$,
$a\da_{A}A$ ($a$ is free from $A$ over $A$).
And since we assume that the class is
$d^{p}$-superstable, it turns out that it is enough to
assume this just for finite $A$, i.e., that the class is weakly simple.

In section~\ref{sec:pre} we give a short summary of the key properties of our context. In section~\ref{sec:sheaf} we illustrate one use of perturbations as a viewpoint to structural approximation. Section~\ref{sec:hom} is a short summary of results from homogeneous model theory that we make use of throughout the paper.

In section~\ref{sec:measuring}
we introduce the main notion of the paper. This is the metric independence notion $a\da^{\e}_{A}B$, $\e >0$,
which intuitively means that the amount of the dependence
of $a$ from $B$ over $A$ is at most $\e$. We start by proving some basic properties
for this notion, e.g., that for all $\e >0$, $a$ and $B$ there is a finite
$A\subseteq B$ such that $a\da_{A}^{\e}B$. As a corollary
we see that for all $a$ and $B$ there is a countable
$A\subseteq B$ such that $a\da_{A}B$. In particular, the class is simple. This is in line with many (discrete)
superstable non-elementary classes in which it also holds that
weak simplicity implies simplicity, see e.g.\ \cite{hykelascartypes}.
However, this is not true if instead of $\distp$-superstability one assumes mere stability, and we give an example of this.

In section~\ref{sec:esplit} we define  Lascar $\e$-splitting and 
study the connection of it with  $\da^{\e}$. Using a Lascar-splitting characterisation we are able to prove a form of monotonicity for $\da^{\e}$ which finally lets us prove that $a\da_AB$ holds if and only if for all $\e>0$, $a\da^{\e}_AB$ holds.

In section~\ref{sec:almost} we investigate a property of the perturbation system which we call almost summability and which improves the behaviour of $\da^{\e}$ over infinite parameter sets. After this we turn our attention to ranks. In section~\ref{sec:entropy} we illustrate how the notion of entropy arises naturally in our context. The way it measures forking chains can be seen as a variant of U-rank. In section~\ref{sec:pregeom} we observe that a certain behaviour of
$\da^{\e}$ in the collection of extensions of
a type,
gives a pregeometry in $\M^{eq}$ (the version used in homogeneous model theory). Here we define another rank function that plays a role in the characterisation.

In the final section~\ref{sec:ex}  we give a worked-through example of the theory
created in this paper. The example is a subclass of Abelian groups
equipped with p-adic metric. It also serves as an example of finding a pregeometry in $\M^{eq}$ although the type we investigate need not be regular.

\section{Assumptions and prerequisites}\label{sec:pre}

Throughout this paper we will  assume $\K$ is a homogeneous metric abstract elementary class with perturbations with complete type spaces that is weakly simple and $d^p$-superstable. We will shortly describe the framework, for full definitions, the reader is referred to \cite{a2}.

The structures studied are many-sorted structures, each sort of which is a complete metric space. One of the sorts is a copy of the ordered field of real numbers. We do not assume the metric spaces to be bounded, nor do we assume uniform continuity of the functions, but the requirements of the class will put some implicit demands on the behaviour of the functions.

We write $a,b$ etc.\ for finite tuples. As a shorthand $ab$ will denote the concatenated tuple of $a$ and $b$. For sets $A$ and $B$, $AB$ will denote their union. We do not specify the sort, but $a\in\M$ will mean that $a$ is a finite tuple of elements of appropriate sorts of $\M$. The distance of two tuples $a$ and $b$ (of the same length) is the maximum of the coordinatewise distances.

A metric abstract elementary class, MAEC, is just a metric adaptation of Shelah's abstract elementary classes, i.e., it consists of a class $\K$ of metric structures and a strong submodel relation $\preccurlyeq$ refining the submodel relation and partially ordering $\K$, such that it satisfies the coherence axiom, the class is closed under (completions) of unions of chains and has a (metric) L\"owenheim-Skolem number (i.e.\ a L\"owenheim-Skolem number where size is measured by density character).

To a given MAEC $(\K,\preccurlyeq)$ we add classes $\BF_\e$ of $\e$-isomorphisms. The $\e$-isomorphisms are uniformly continuous bijections between members of $\K$ satisfying natural properties of composition and inversion, with $\e$ being a measure for the amount of error the functions in $\BF_\e$ make. A natural example of an $\e$-isomorphism is a linear isomorphism $T$ of a Banach space with $\norm{T},\norm{T^-1}\leq e^\e$.

We further assume that the class has arbitrarily large models, satisfies joint embedding (JEP), amalgamation wrt.\ the $\e$-isomorphisms (AP) and is homogeneous. This allows us, for any given cardinal $\mu$, to construct a $\mu$-universal, strongly $\mu$-homogeneous monster model $\Mon$ with the additional property, that any $\e$-isomorphism between strong submodels of $\M$ extends to an $\e$-automorphism of $\M$. The notion of type used is Galois-type, $\tg(a/A)$, which corresponds to the orbit of $a$ under automorphisms of $\M$ fixing $A$ pointwise.

To define a metric on the type space we define the relation  $\distp$:
\begin{definition}\label{dist}
For tuples $a,b\in\M$ and $\varepsilon>0$ we write
$$
\distp(\tg(a/\emptyset),\tg(b/\emptyset))\leq\varepsilon
$$ 
if there are 
$\varepsilon$-automorphisms $f$ and $g$ of $\M$ such that
$\de(f(a),b)\leq\varepsilon$ and $\de(g(b),a)\leq\varepsilon$.

For types over parameter sets $A$, we define
$$
\distp(\tg(a/A),\tg(b/A))=\sup\{\distp(\tg(ac/\emptyset),\tg(bc/\emptyset)):
c\in A \textrm{ finite}\}.
$$
\end{definition}

To get a metric, we assume the following perturbation property:  whenever $a,b\in\Mon$
are tuples such that $\distp(\tg(a/\emptyset),\tg(b/\emptyset))=0$
(i.e. $\distp(\tg(a/\emptyset),\tg(b/\emptyset))\leq\varepsilon$ for all
positive $\varepsilon$) then $\tg(a/\emptyset)=\tg(b/\emptyset)$.

Strictly speaking, $\distp$ is not a metric (the triangle inequality needs some rescaling), but it defines a metrisable uniformity, so it makes sense to talk about Cauchy sequences, limits and completeness with respect to $\distp$.

Finally, the assumption of \emph{complete type-spaces} states, that all $\distp$-Cauchy sequences over $\emptyset$ have a limit. A fact that was not pointed out in \cite{a2} is:

\begin{fact}\label{fact6}
If $\K$ is homogeneous, satisfies the perturbation property and has complete type spaces then $d^{p}$-Cauchy sequences over any parameter set converge.
\end{fact}

Stability is then defined in the natural way:
$\K$ is \emph{$\lambda$-$\distp$-stable} if for any
set $A$ with $|A|\leq\lambda$, the set of types over $A$ has density
$\lambda$ with respect to
$\distp$.

Note that the cardinality of a type space of density $\lambda$ (wrt.\ $\distp$) has cardinality $\leq\lambda^\omega$, so the notion of \emph{stability} is the same regardless of whether we consider the class as a (perturbed) metric class or as a discrete AEC, although the cardinalities where stability occurs vary.

\section{An example of the use of perturbations}\label{sec:sheaf}

This section describes an example of the use of $\e$-isomorphisms in model theory. It arises from the sheaf of
rational Weyl algebras studied in \cite{zilber:Walg}
by B. Zilber. In the end of the example we show how
one can get a natural sheaf topology 
from the perturbations of the class.
In this example we leave the details to the reader to check.
In our final example
in the end of this paper, all the details are given.

Weyl algebras are algebras generated by $P$ and $Q$
(representing momentum and position in quantum mechanics),
with the canonical commutation relation $QP-PQ=i\hbar$.
This cannot be represented by bounded operators in a Hilbert space,
so instead operators $U^t=\exp(itQ)$ and $V^v=\exp(ivP)$ are studied.
Zilber studies the algebras arising from
letting the commutator of $U$ and $V$ be a root of unity,
choosing $t$ and $v$ to
be
rational
and by these approximating the full Weyl algebra. We show how
perturbations give a viewpoint to this idea of structural approximation.

The idea is to build a large Hilbert space over an indexing set $R$. To each $a\in R$ we attach an orthonormal sequence of eigenvectors of $U$ such that $V$ acts as a shifting operator on (all or part of) this sequence.

We let $\K$ be the class of structures $\A$ that consists of
the following two parts:
\begin{itemize}
\item[(a)] A set $R$ with discrete metric (distance between any two
points is $1$).
\item[(b)] A complex Hilbert space $H$ with the usual metric that comes from
the inner product.
\end{itemize}
$R$ works as an indexing set. In order to put these into the many-sorted form demanded, we actually include four sorts: $R, H,\C$ and $\R$.

In addition to these, in $\A$,
there are interpretations for the unary
function symbols: $U$, $V$ and $B_{z}$, for all rational numbers $z$:
The functions $B_{z}$ are from $R$ to $H$
and the requirement is
that $\{ B_{z}(a)\vert\ a\in R,\ z\in\Q\}$
forms an orthonormal basis of $H$.
The functions $U$ and $V$ are unitary operators on $H$ but to
give the requirements we attach
numbers to the elements of $R$. The numbers are not part of
the structure but in the end they can be read out from it.
So to each $a\in R$ we attach four numbers: $N_{a}$, which is a
positive natural number such that for all
$a\in R$, $N_a=n!$ for some positive natural number $n$
(the latter requirement is just to make fractions easier). 
We have also
a complex number
$q_{a}$ of absolute value $1$. Finally we have
natural numbers $n_{a}$ and $m_{a}$.
If $q_{a}$ is not a root of unity or is $1$,
then $n_{a}=m_{a}=0$. Otherwise
$q_{a}=e^{i2\pi m/n}$ for some natural numbers $0<m<n$ and 
then $m_{a}/n_{a}=m/N_a^{2}n$ and $gcd(m_{a},n_{a})=1$.
To simplify the notation, for $q=e^{i2\pi r}$, $0\le r<1$,
and integer numbers $n\ne 0$ and $m$,
we write $q^{m/n}$ for $e^{i2\pi rm/n}$.

Now we are ready to state the requirements for
the operators $U$ and $V$.
If $n_{a}=0$, then for all integers $z$,
$U(B_{z/N_{a}^{2}}(a))=q_ {a}^{z/N_{a}^{2}}B_{z/N_{a}^{2}}(a)$ and
$V(B_{z/N_{a}^{2}}(a))=B_{(z-1)/N_{a}^{2}}(a)$, and
for all other rational numbers $r$, $U(B_{r}(a))=V(B_{r}(a))=B_{r}(a)$.
If $n_{a}\ne 0$, let $z_{a}$ and $z'_{a}$ be integers such that
$z_{a}+z'_{a}\in\{ 0,1\}$ and $z'_{a}-z_{a}+1=n_{a}$.
Then if $z_{a}\le z\le z'_{a}$, $U(B_{z/N_{a}^{2}}(a))=
q_ {a}^{z/N_{a}^{2}}B_{z/N_{a}^{2}}(a)$ 
and for all other rational numbers $r$
$U(B_{r}(a))=B_{r}(a)$.
If $z_{a}<z\le z'_{a}$, then $V(B_{z/N_{a}^{2}}(a))=B_{(z-1)/N_{a}^{2}}(a)$,
$V(B_{z_{a}/N_{a}^{2}}(a))=B_{z'_{a}/N_{a}^{2}}(a)$ and
for all other rational numbers $r$, $V(B_{r}(a))=B_{r}(a)$.

If as a strong submodel relation we use the submodel relation,
$\K$ is easily seen to be a homogeneous MAEC with AP and JEP,
the first-order theory of the monster model is unstable
and as a MAEC the class is superstable but not $\o$-stable.
Notice that if we write
$H_{a}$ for the subspace of $H$ generated by
$\{ B_{r}(a)\vert\ r\in\Q ,\ V(B_{r}(a))\ne B_{r}(a)\}$,
then $U$ and $V$ restricted to this space satisfy
the Weyl commutator law $V^{z}U^{z'}=q_{a}^{zz'/N_{a}^{2}}U^{z'}V^{z}$
for all integers $z$ and $z'$ (so as a representation of a Weyl algebra,
$U$ and $V$ restricted to $H_{a}$ should really be understood as
$U^{1/N_{a}}$ and $V^{1/N_{a}}$).

Let us then define the sets $\BF_{\e}$ of generalized isomorphisms:
We put $f:\A\rightarrow\B$ in $\BF_{\e}$
if $f$ is an isomorphism after we remove the interpretations
of $U$ and $V$ from the structures and in addition
the following holds:
For all $a\in R^{\A}$ and $b\in R^{\B}$, if $b=f(a)$, then
$$(*)\ \ \ \vert 1/N_{a}-1/N_{b}\vert
+\vert q_{a}-q_{b}\vert +\vert 1/n_{a}-1/n_{b}\vert
\le\e ,$$
where $1/0$ is considered to be $0$.
It is easy to see that $(\K ,\subseteq ,\BF_{\e})_{\e\ge 0}$
is a homogeneous MAEC with perturbations, it has JEP and AP (for
the perturbations) and the perturbation property. With these perturbations,
the class is $\o$-stable.

However, $\K$ does not have a complete type space:
Let $\M$ be the monster model of this class and choose
$a_{i}\in R$, $i<\o$, so that for all $i<j<\o$,
$q_{a_{i}}=q_{a_{j}}$ and $N_{a_{i}}=(i+1)!$.
Let $p_{i}=t^{g}(a_{i}/\empty )$.
Then the types $p_{i}$ form a Cauchy sequence without a limit.
Of course, we can get a complete type space by replacing (*) above by
$$(**)\ \ \ \vert N_{a}-N_{b}\vert
+\vert q_{a}-q_{b}\vert +\vert 1/n_{a}-1/n_{b}\vert
\le\e.$$
We did not do it this way because now we can discuss
the limit problem: Another natural way of fixing this problem
of complete type space would
be by adding the limits to the structure.
In his sheafs, Zilber studies limits of this kind
(Zilber looks only at the cases in which the commutator is a root of unity
and $\ne 1$).
We choose not to add any such limits here
because we see problems in choosing these
limits from the point of view of physics (we study this in \cite{qm},
and if one wants to add these limit one needs to replace $U$ and $V$
with the operators $U^{1/n}$ and $V^{1/n}$ for all positive natural numbers
$n$).

Now suppose $a,a_{i}\in R$ are such that for all $i<\o$,
there is $f_{i}\in \BF_{1/(i+1)}$ mapping $a_{i}$ to $a$ and
for all $i<j<\o$, $q_{a_{i}}\ne q_{a_{j}}$.
Notice that then $n_{a}=0$ and $N_{a_{i}}=N_{a}$ for all large enough
$i$. Also if we write $H^{*}_{a}$ for the subspace of $H$
generated by $\{ B_{r}(a)\vert\ r\in\Q\}$ and $U_{a}$ and $V_{a}$
for the restrictions of $U$ and $V$ to $H^{*}_{a}$,
then $f_{i}U_{a_{i}}f_{i}^{-1}$ converges to $U_{a}$
and $f_{i}V_{a_{i}}f_{i}^{-1}$ converges to $V_{a}$
in the weak topology (not in the operator norm).
If one defines the time evolution operators to these
spaces so that the same happens with the time evolution operator,
this convergence is strong enough to allow one to calculate
e.g.\ propagators in the case the Weyl commutator is not a root
of unity by calculating them in the root of unity cases,
i.e., in finite dimensional cases in which number theory
can be used to do the calculations and in which we have `all'
eigenvectors. Notice that if $q_{a}$ is not a root of unity,
then $V\raj H_{a}$ does not have any eigenvectors,
and thus neither does the time evolution operator (at least in most cases),
which makes the direct calculations very hard.

We finish this example by showing how to get a natural sheaf
from the perturbations of our class. So we work inside e.g.\ the
monster model $\M$ of the class $\K$.

As the basis of the sheaf, we choose $R$ and as a topology
on $R$ we use the one
we get from the pseudometric $d$, where $d(a,b)$ is the infimum
of all $\e\ge 0$ such that there is an $\e$-automorphism
$f$ of $\M$ such that
$f(a)=b$. The set of elements of the sheaf is $E=\bigcup_{a\in R}H^{*}_{a}$
and the topology on $E$
is given by the pseudometric  $d^{*}$
where $d^{*}(x,y)=10$ if there is no $\e$-automorphism $f$
such that $f(x)=y$
for any $\e\ge 0$, and otherwise $d^{*}(x,y)$ is the infimum
of all $\e\ge 0$ such that there is an $\e$-automorphism
$f$ such that
$f(x)=y$. The projection $p:E\rightarrow R$ is the obvious one:
$p(x)=a$ if $x\in H^{*}_{a}$. Clearly $p$ is a local homeomorphism.
In \cite{ovsheaf}, model theory for sheafs of metric structures is studied
and our sheaf $E$ seems to fit into their framework.

\section{Using discrete homogeneous model theory results}\label{sec:hom}

In this section we
note  that by Morleyisation (introducing predicates for the types) we may turn $\Mon$ into a homogeneous first order structure, preserving stability, and thus use results from \cite{hysh}.

By $\l (\K )$ we mean the least cardinal $\l$ such that
as a homogeneous AEC, $\K$ is $\l$-stable (i.e. $\l$-stable
in the sense of \cite{hysh}).
By $\kappa(\K)$ we denote the least cardinal $\kappa$ such that there is no strongly splitting sequence of length $\kappa$.

We say that $a$ and $b$ have the same Lascar strong types over $A$, $Lstp(a/A)=Lstp(b/A)$, if  $E(a,b)$ holds for any $A$-invariant equivalence relation with a bounded number of equivalence classes.

\begin{fact}\label{fact1}
Let $\M$ be strongly $\lambda$-homogeneous. For every $\kappa<\lambda$ there is a cardinal $H(\kappa)$ 
such that if $A$ is a set of size $\leq\kappa$ and $(a_i)_{i<H(\kappa)}\subset\M$ then there exists an $A$-indiscernible sequence $(b_i)_{i<\omega}\in\M$ such that for every $n<\omega$ there exist $i_0<\cdots<i_n<H(\kappa)$ such that
$$t^g(b_0,\dots,b_n/A)=t^g(a_{i_0},\dots,a_{i_n}/A).$$
\end{fact}

Note that the fact ensures that over any set $A$ there are less than $H(\abs{A})$ Lascar strong types over $A$.

In a stable homogeneous class we can define an independence notion based on strong splitting as done in \cite{hysh}: 

We write $a\da_AB$ if there is $C\subseteq A$ of power $<\kappa(\K)$ such that for all $D\supseteq A\cup B$ there is $b$ which satisfies $t^g(b/AB)=t^g(a/AB)$ such that $t^g(b/D)$ does not split strongly over $C$. For an arbitrary set $C$, $C\da_AB$ means $a\da_AB$ for all finite tuples $a\in C$.

\begin{fact}[Hyttinen-Shelah \cite{hysh}]\label{fact2}
In a stable homogeneous class $\da$ satisfies:
\begin{itemize}
\item[(i)] (monotonicity) If $A\subseteq A'\subseteq B'\subseteq B$ and $a\da_AB$ then $a\da_{A'}B'$.
\item[(ii)] (extension of free types) If $A\subseteq B$, $a\da_AA$ and $t^g(a/A)$ is unbounded, then there is $b$ such that $b\da_AB$ and $Lstp(b/A)=Lstp(a/A)$.
\item[(iii)] (finite character) If $A\subseteq B$, $a\nda_AB$ and $a\da_AA$ and $t^g(a/A)$ is unbounded then there is some finite $B'\subseteq B$ such that $a\nda_AB'$.
\item[(iv)] (symmetry for free types) For all $a$, $b$ and $A$, $b\da_AA$ and $a\da_Ab$ implies $b\da_Aa$. By finite character this generalises to: if $A\da_BC$ and $C\da_BB$ then $C\da_BA$.
\item[(v)] (pair) If $b\da_AD$ and $c\da_{Ab}D$ then $bc\da_AD$.
\item[(vi)] (stationarity) If $a\da_Ac$, $b\da_Ac$ and $Lstp(a/A)=Lstp(b/A)$ then $t^g(a/Ac)=t^g(b/Ac)$.
\end{itemize}
\end{fact}

\begin{lemma}\label{lemma3}
If $C\da_AB$ and $D\da_{AC}B$ then $CD\da_AB$.
\end{lemma}

\begin{proof}
By strong extension and stationarity of Lascar strong types we may assume $B$ is $\lambda(\K)$-saturated. 
We may also assume $A$ and $C$ are of power $<\kappa(\K)$.
Now if $CD\nda_AB$ there are finite $c\in C$ and $d\in D$ such that $cd\nda_AB$. So there is an $A$-indiscernible sequence $I=(a_i)_{i<\omega}\subset B$ such that $t^g(cda_0/A)\neq t^g(cda_1/A)$. If $I$ is $AC$-indiscernible, this contradicts $D\da_{AC}B$. So $I$ cannot be indiscernible over $AC$ but then (by re-enumerating) for some $n$ and $c\in C$ $t^g(c,a_0,\dots,a_{n-1}/A)\neq t^g(c,a_n,\dots,a_{2n-1}/A)$ giving an $A$-indiscernible sequence contradicting $C\da_AB$.
\end{proof}

\begin{corollary}\label{corollary4}
If $A\subseteq B$, $a\da_AB$, $a\da_BC$, $C\da_BB$ and $B\da_AA$ then $a\da_AC$.
\end{corollary}

\begin{proof}
 By symmetry $B\da_Aa$ and $C\da_Ba$ and thus by Lemma~\ref{lemma3}
$BC\da_Aa$. By symmetry we then have $a\da_ABC$. 
\end{proof}

In this context we define weak simplicity as $a\da_AA$ for all $a$ and finite $A$. Thus if $\K$ is stable and weakly simple then $\da$ satisfies monotonicity and stationarity of strong types, and over finite sets in addition transitivity, symmetry and strong extension.

As in \cite{hykesuperstability}, we write $Lstp^{w}(a/A)=Lstp^{w}(b/A)$
if for all finite $B\subseteq A$,
$Lstp(a/B)=Lstp(b/B)$. Following \cite{hykelascartypes}, the types
$Lstp^{w}$ are called Lascar types.
By homogeneity, 
if $Lstp^{w}(a/A)=Lstp^{w}(b/A)$ then
$t^g(a/A)=t^g(b/A)$.

The following lemma is needed because the $d^{p}$-distance of
Galois-types we use need not be a metric, see \cite{a2}.
If there are no perturbations, i.e.
$F_{\e}=F_{0}$ for all $\e >0$, then it is and we can
choose $n(\e )=\e /n$ and $\d_n=2^{-n}$.

\begin{lemma}\label{lemma5}
\begin{itemize}
\item[(i)] For all $n>1$ and $\e >0$, there is $n(\e )>0$ such that
for all $a_{i}$, $i\le n$, if for all $i<n$
$d^{p}(t^g(a_{i}/\empty ),t^g(a_{i+1}/\empty ))\le n(\e )$,
then $d^{p}(t^g(a_{0}/\empty ),t^g(a_{n}/\empty ))\le\e$.
\item[(ii)] There are $\d_n>0$ such that if $d^p(t^g(a_n/\empty),t^g(a_{n+1}/\empty))\le\d_n$ for all $n<\o$ then $(t^g(a_i/A))$ is a Cauchy sequence (wrt.\ $d^p$).
\end{itemize}
\end{lemma}

\begin{proof}
 Immediate by the definitions, see \cite{a2}. 
\end{proof}

\section{Measuring independence}\label{sec:measuring}

In this section we define a distance-like relation $d^p_a$ on the space of Lascar types. As for $d^p$ in \cite{a2} it is not exactly a metric but defines a metrisable topology. Using $d^p_a$ we define $\e$-independence (Definition~\ref{definition10}) and explore its properties.

Recall that throughout the paper we assume $\K$ is a homogeneous MAEC with perturbations with complete type spaces that is weakly simple and $d^p$-superstable (except in Corollary~\ref{corollary14} where we give a characterisation of superstability, and thus only assume stability). Note, however, that completeness of type-spaces is used only to prove extension in Lemma~\ref{lemma19}, so we could omit the assumption if we instead choose to assume extension.

The aim is to show that $\da^{0}$ and $\da$ agree and that the class actually is simple, but to get there we need many intermediate results.
Many of the classical properties of independence (monotonicity, transitivity), do not hold for $\e$-independence, but for $0$-independence, and to show this we look at the interplay between $\e$-independence and ordinary independence (as defined for homogeneous discrete classes). Until we prove simplicity (in Corollary~\ref{corollary23}) we often have to work over finite sets, as these are the only ones over which weak simplicity guarantees that $\da$ is well behaved.

\begin{definition}\label{definition7}
\begin{itemize}
\item[(i)] For a finite set $A$ we define
$$d^{p}_{a}(Lstp(a/A),Lstp(b/A))=\sup\{d^{p}(t^g(a/B),t^g(b/B)):A\subseteq B \, \textrm{finite},\, B\da_Aab\}.$$
\item[(ii)] For any set $B$ we then define
$$
d^{p}_{a}(Lstp^w(a/B),Lstp^w(b/B)) = \sup\{d^{p}_{a}(Lstp(a/A),Lstp(b/A)):A\subseteq B,\, A\, \textrm{finite}\}.$$
\end{itemize}
\end{definition}

We first show that this definition makes sense:

\begin{lemma}\label{lemma8}
\begin{itemize}
\item[(i)] If $A$ is finite, then the definitions (i) and (ii) of \ref{definition7} give the same result.
\item[(ii)] If $Lstp^{w}(a/B)=Lstp^{w}(a'/B)$
and $Lstp^{w}(b/B)=Lstp^{w}(b'/B)$, then
$$d^{p}_{a}(Lstp^{w}(a/B),Lstp^{w}(b/B))=
d^{p}_{a}(Lstp^{w}(a'/B),Lstp^{w}(b'/B)).$$
\item[(iii)] If $A$ is finite and $ab\da_AB$ then 
$$
d^p_a(Lstp^w(a/B),Lstp^w(b/B))= d^p_a(Lstp(a/A),Lstp(b/A)).
$$
\item[(iv)] If $d^{p}_{a}(Lstp^{w}(a/B),Lstp^{w}(b/B))=0$,
then $Lstp^{w}(a/B)=Lstp^{w}(b/B)$.
\end{itemize}
\end{lemma}

\begin{proof}
(i)--(iii): Immediate by the definitions.

(iv): It suffices to show that for any finite $A$,
if $d^{p}_{a}(Lstp(a/A),Lstp(b/A))=0$,
then $Lstp(a/A)=Lstp(b/A)$. For this choose
$c$ so that $Lstp(c/A)=Lstp(a/A)$ and
$c\da_{A}ab$. Then by the assumption,
$d^{p}(t^g(a/Ac),t^g(b/Ac))=0$ and thus by the perturbation
property, $t^g(a/Ac)=t^g(b/Ac)$ and thus
$Lstp(a/A)=Lstp(b/A)$. 
\end{proof}

Although $d^{p}_{a}$ may not satisfy the triangle
inequality, as in \cite{a2} for $d^{p}$,
it gives rise to a metrisable topology to the set of
all Lascar types over any fixed set $B$. In fact we have the following analogue of Lemma~\ref{lemma5} (i):

\begin{lemma}\label{lemma9}
For all $n>1$ and $\e >0$, for all $a_{i}$, $i\le n$, and all $A$ if for all $i<n$ $d^{p}_{a}(Lstp^w(a_{i}/A ),Lstp^w(a_{i+1}/A))\le n(\e )$, where $n(\e)$ is as in Lemma~\ref{lemma5} (i), 
then $d^{p}_{a}(Lstp^w(a_{0}/A ),Lstp^w(a_{n}/A))\le\e$.
\end{lemma}

\begin{proof}
 It suffices to prove this when $A$ is finite. For this assume that $d^{p}_{a}(Lstp(a_{i}/A ),Lstp(a_{i+1}/A))\le n(\e )$ and let $D\supseteq A$ be finite and such that $D\da_Aa_0a_n$. Choose $D'\supseteq A$ satisfying $Lstp(D'/Aa_0a_n)=Lstp(D/Aa_0a_n)$ and $D'\da_A\bigcup_{i\leq n}a_i$. Then by assumption $d^p(t^g(a_i/D'),t^g(a_{i+1}/D'))\leq n(\e)$ and thus by Lemma~\ref{lemma5} (i), $d^p(t^g(a_0/D'),t^g(a_n/D'))=d^p(t^g(a_0D'/\empty),t^g(a_nD'/\empty))\leq\e$. As $t^g(D'/Aa_0a_n)=t^g(D/Aa_0a_n)$ we are done. 
\end{proof}

We are ready to define the main notion of this paper:

\begin{definition}\label{definition10}
For $\e >0$, we write
$a\da^{\e}_{A}B$ if for all finite $C\subseteq A$,
there is some finite $D$ with $C\subseteq D\subseteq A$ and
$b$ such that
$Lstp(b/D)=Lstp(a/D)$, $b\da_{D}AB$ and
$d^{p}_{a}(Lstp^{w}(b/AB),Lstp^{w}(a/AB))\le\e$.
By $a\da^{0}_{A}B$ we mean that $a\da^{\e}_{A}B$ holds
for all $\e >0$.
\end{definition}

This independence notion has some immediate properties. Note, however, that we only have partial monotonicity.

\begin{lemma}\label{lemma11}
Suppose $\e >0$.
\begin{itemize}
\item[(i)] If $A\subseteq C\subseteq D$
and $a\da^{\e}_{A}D$, then $a\da^{\e}_{A}C$.
\item[(ii)] If $ab\da^{\e}_{A}B$, then $a\da^{\e}_{A}B$.
\item[(iii)] If $a\da^{\e}_{A}B$, $A'\subseteq A$ is finite, then there is some finite $A''$ such that $A'\subseteq A''\subseteq A$ and $a\da^{\e}_{A''}B$.
\item[(iv)] If $A$ is finite and
$a\nda^{\e}_{A}B$, then there is a finite $C\subseteq B$
such that $a\nda^{\e}_{A}C$.
\item[(v)] If $\e >\d >0$, then
$a\da^{\d}_{A}B$ implies $a\da^{\e}_{A}B$.
\item[(vi)] If $A$ is finite and $a\da_{A}B$, then $a\da^{0}_{A}B$.
\item[(vii)] If $A$ is finite, $a\da_AB$ and $a\da^{\e}_{AB}C$ then $a\da^{\e}_ABC$.
\item[(viii)] If $A$ is finite, $a\da_AB$ and $a\da^{\e}_ABC$ then $a\da^{\e}_{AB}C$.
\end{itemize}
\end{lemma}

\begin{proof}
Immediate by the definitions.
\end{proof}

Now $d^p$-superstability ensures that $\da^{\e}$ has local character:

\begin{lemma}\label{lemma12}
\begin{itemize}
\item[(i)] For all $\e >0$, $a$ and $A$ there is a finite $B\subseteq A$
such that $a\da^{\e}_{B}A$.
\item[(ii)] For all $a$ and $A$, $a\da^{0}_{A}A$.
\end{itemize}
\end{lemma}

\begin{proof}
(i): By Lemma~\ref{lemma11}, it is enough to show that
there are no $a$ and finite $A_{i}$, $i<\o$,
such that for all $i<\o$, $A_{i}\subseteq A_{i+1}$
and $a\nda^{\e}_{A_{i}}A_{i+1}$. For a contradiction,
suppose such $a$ and $A_{i}$, $i<\o$, exist.

Let $\k >\l (\K )$ be a cardinal of cofinality $\o$
and such that $\K$ is $\k$-$d^{p}$-stable. 

We define a new increasing sequence of finite sets $A'_i$ such
that every $b\models t^g(a/A'_i)$ with
$b\da_{A'_i}A'_{i+1}$ satisfies $d^p(t^g(b/A'_{i+1}),t^g(a/A'_{i+1}))>\e$:
For each $i<\o$ let $c_i\models Lstp(a/A_i)$ with $c_i\da_{A_i}a$.
Now let $A'_0=A_0c_0$. When $A'_i$ has been defined
such that $b\models t^g(a/A'_i)$ implies $b\models Lstp(a/A_i)$
and $a\da_{A_i}A'_i$ let $b_i\models Lstp(a/A_i)$, $b_i\da_{A_i}A_{i+1}$.
As $a\nda^{\e}_{A_i}A_{i+1}$, we have $d^p_a(Lstp(b_i/A_{i+1}),Lstp(a/A_{i+1}))>
\e$, i.e., there is some finite $B_i\supseteq A_{i+1}$ with $B_i\da_{A_{i+1}}ab_i$ such that $d^p(t^g(a/B_i),t^g(b_i/B_i))>\e$ and we may assume $B_i\da_{A_{i+1}}ab_ic_{i+1}$ 
. Then define $A'_{i+1}=A'_iB_ic_{i+1}$. 

Now if $b\models t^g(a/A'_i)$ and $b\da_{A'_i}A'_{i+1}$, these imply $b\models Lstp(a/A_i)$ and $b\da_{A_i}B_i$. So $b\models t^g(b_i/B_i)$ and thus $d^p(t^g(a/B_i),t^g(b/B_i))>\e$.

Then we can proceed with the usual construction from \cite{she}: 
For all $\n\in\k^\o$ and all $n<\o$, choose $A_{\n\raj n}$ and $a_\n$ so that
\begin{itemize}
\item[(a)] for all $\n\in \k^{\o}$ there is an automorphism $F_{\n}$
of the monster model such that $F_{\n}(a_{\n})=a$,
for all $i<\o$, $F_{\n}(A_{\n\raj i})=A'_{i}$ and if
$\xi\in \k^{\o}$ and
$\n\raj i=\xi\raj i$, then $F_{\n}\raj A_{\n\raj i}=
F_{\xi}\raj A_{\n\raj i}$,
\item[(b)] for all $\n\in \k^{\o}$ and $i<\o$,
$$a_{\n}\da_{A_{\n\raj i}}\cup\{A_{\xi}\vert\ \xi\in \k^{<\o},
\ \n\raj i+1\not\subseteq\xi\}.$$
\end{itemize}
Let $D =\cup_{\n\in \k^{<\o}}A_{\n}$. Now
clearly $d^{p}(t^g(a_{\n}/D ),t^g(a_{\xi}/D ))>\e$
for distinct $\n ,\xi\in \k^{\o}$. This contradicts
the choice of $\k$.

(ii): Note that if $a\nda^{\e}_{A}A$ then there is a finite $A'\subseteq A$ such that $a\nda^{\e}_{B}A$ for all finite $B$ with  $A'\subseteq B\subseteq A$. Then proceed as in (i). 
\end{proof}

\begin{corollary}\label{corollary13}
For all $A$ and $a$, there is countable $B\subseteq A$ such that $a\da^{0}_BA$.
\end{corollary}

\begin{proof}
 We define $B$ by induction. Let $B_0=\emptyset$.  When a finite $B_n$ has been defined we define $B_{n+1}\supseteq B_n$ finite such that $a\da^{1/(n+1)}_{B_{n+1}}A$: By Lemma~\ref{lemma12} (ii) $A\da^{0}_AA$ and thus by Lemma~\ref{lemma11} (iii) there is some finite $B_{n+1}$ with $B_n\subseteq B_{n+1}\subseteq A$ such that $a\da^{1/(n+1)}_{B_{n+1}}A$. In the end let $B=\bigcup_{n<\omega}B_n$. Then for any $\e>0$ and finite $C\subseteq B$ there is $n>1/\e$ such that $C\subseteq B_n$ and $B_n$ witnesses (as $D$ in Definition~\ref{definition10}) $a\da^{\e}_BA$. 
\end{proof}

Now we can see that $d^p$-superstability has natural characterisations in terms of  $\da^{\e}$-forking sequences:

\begin{corollary}\label{corollary14}
Suppose the class $\K$ is stable and weakly simple. Then T.F.A.E.
\begin{itemize}
\item[(i)] $\K$ is $d^p$-superstable.
\item[(ii)] For no $\e>0$ is there an infinite $\da^{\e}$-forking sequence.
\item[(iii)] For all $a$, $A$ and $\e>0$, there is a finite $B\subseteq A$ such that $a\da^{\e}_BA$.
\end{itemize}
\end{corollary}

\begin{proof}
 (i)$\Rightarrow$(ii)$\Rightarrow$(iii) follow by Lemmas \ref{lemma12} and \ref{lemma11}
so we prove (iii)$\Rightarrow$(i).

Let $H(\aleph_0)$ be as in Fact~\ref{fact1}. We claim that $\K$ is $d^p$-stable in every
$\kappa\geq H(\aleph_0)$. In fact, the density character of the set of Lascar types ($Lstp^w$) over a set $A$ (wrt.\ $d^p_a$) is at most $\abs{A}+H(\aleph_0)$.

So let $\abs{A}=\kappa\geq H(\aleph_0)$. If the density character of the set
of Lascar types over $A$ is greater than $\kappa$ there are some $\e>0$
and tuples $a_i$ for $i<\kappa^+$ such that
$d^p_a(Lstp^w(a_i/A),Lstp^w(a_j/A))>\e$ for all $i\neq j$.
Let $\delta=2(\e)$ (from Lemma~\ref{lemma9}). By (iii) there are finite sets $A_i\subset A$ such
that $a_i\da^{\delta}_{A_i}A$. Since there are only $\kappa$ finite
subsets of $A$, $\kappa^+$ many of the $A_i$s are the same set $A'$.
Further, since there are less than $H(\aleph_0)$ Lascar 
types over $A'$, for $\kappa^+$ many indices the Lascar strong type
$Lstp(a_i/A')$ is the same. Let $Lstp(a/A')=Lstp(a_i/A')$ and $a\da_{A'}A$.
Then for $\kappa^+$ many indices $d^p_a(Lstp^w(a/A),Lstp^w(a_i/A))\leq\delta$
and thus by Lemma~\ref{lemma9}
$d^p_a(Lstp^w(a_i/A),Lstp^w(a_j/A))\leq\e$ for $\kappa^+$ many $i,j<\kappa^+$, a contradiction. 
\end{proof}

Next we have a look at versions of transitivity.

\begin{lemma}\label{lemma15}
Suppose $A$ and $B$ are finite, $\e >0$, $a\da^{\e}_{A}B$ and $a\da_{AB}C$.
Then $a\da^{\e}_{A}BC$.
\end{lemma}

\begin{proof}
Clearly it is enough to prove this for such $A$, $B$ and $C$ that $A\subseteq B\subseteq C$ and $C$ is finite. For this let $b$ be such that
$Lstp(b/A)=Lstp(a/A)$ and $b\da_{A}C$.
We need to prove that
$d^{p}_{a}(Lstp(b/C),Lstp(a/C))\le\e$. By Lemma~\ref{lemma8}(ii) 
$d^{p}_{a}(Lstp(c/C),Lstp(a/C))=d^{p}_{a}(Lstp(b/C),Lstp(a/C))$
for all $c$ such that $Lstp(c/C)=Lstp(b/C)$
and thus we may assume that
$b\da_{A}Ca$. Let $D\supseteq C$ be finite
and such that $D\da_{C}ab$.
We need to show that $d^{p}(t^g(a/D),t^g(b/D))\le\e$.

But now $a\da^{\e}_{A}B$ and thus
$d^{p}_{a}(Lstp(b/B),Lstp(a/B))\le\e$.
By transitivity
$D\da_{B}ab$,  so by Lemma~\ref{lemma8}(iii)  $d^{p}(t^g(a/D),t^g(b/D))\le\e$. 
\end{proof}

\begin{corollary}\label{corollary16}
If $A$ and $B$ are finite then for all $\e >0$, if $a\da^{\e}_{A}B$,
then for all $C$ there is $b$ such that
$b\da^{\e}_{A}BC$ and $Lstp(b/AB)=Lstp(a/AB)$.
\end{corollary}

\begin{proof}
Just choose $b\models Lstp(a/AB)$ satisfying $b\da_{AB}C$ and use Lemma~\ref{lemma15}. 
\end{proof}

\begin{lemma}\label{lemma17}
For all $\e>0$ there exists a $\d>0$ (namely $\d=2(\e)$) such that if $A\subseteq B\subseteq C$, $a\da^{\d}_AB$ and $a\da^{\d}_BC$, then $a\da^{\e}_AC$. In particular $\da^{0}$ satisfies transitivity.
\end{lemma}

\begin{proof}
 Let $\d=2(\e)$. By Lemma~\ref{lemma11} we may assume $A\subseteq B\subseteq C$ are finite. Let $c$ be such that $Lstp(c/A)=Lstp(a/A)$ and $c\da_AC$. We wish to show that $d^p_a(Lstp(a/C),Lstp(c/C))\leq\e$.

Choose $b$ such that $Lstp(b/B)=Lstp(a/B)$ and $b\da_BC$. As $a\da^{\d}_BC$, we must have $d^p_a(Lstp(a/C),Lstp(b/C))\leq\d$. By Lemma~\ref{lemma15}, $b\da^{\d}_AC$, and thus $d^p_a(Lstp(b/C),Lstp(c/C))\leq\d$. But by Lemma~\ref{lemma9},
$d^p_a(Lstp(a/C),Lstp(c/C))\leq \e$ proving $a\da^{\e}_AC$. 
\end{proof}

We can now prove stationarity for Lascar types.

\begin{lemma}\label{lemma18}
Suppose $Lstp^w(a/A)=Lstp^w(b/A)$, $a\da^{0}_AB$ and $b\da^{0}_AB$. Then $Lstp^w(a/B)=Lstp^w(b/B)$.
\end{lemma}

\begin{proof}
Let $\e>0$. W.l.o.g.\ we may assume $A\subseteq B$. We show that $d^p_a(Lstp^w(a/B),Lstp(b/B))\leq\e$. Let $\d\le2(\e)$ and $\d'\leq 2(\d)$. By Lemma~\ref{lemma12} let $C\subseteq A$ be finite such that $ab\da^{\d'}_CA$. Then $a\da^{\d'}_CA$ and $a\da^{0}_AB$ so by Lemma~\ref{lemma17},
$a\da^{\d}_CB$. Similarly $b\da^{\d}_CB$. By assumption $Lstp(a/C)=Lstp(b/C)$ and if we choose $c$ such that $Lstp(c/C)=Lstp(a/C)$ and $c\da_CB$, we have $d^p_a(Lstp^w(a/B),Lstp^w(c/B))\leq\d$ and $d^p_a(Lstp^w(b/B),Lstp^w(c/B))\leq\d$. By Lemma~\ref{lemma9}, 
$d^p_a(Lstp^w(a/B),Lstp^w(b/B))\leq\e$. 
\end{proof}

Next we show extension over countable sets. This lemma is the only place where we use completeness of type spaces, so if we instead wish to assume countable extension we could omit that assumption.

\begin{lemma}\label{lemma19}
For any $a$ and $B$ and any countable $A\subseteq B$ there exists some $a'$ satisfying  $Lstp^w(a'/A)=Lstp^w(a/A)$ and $a'\da^{0}_{A}B$. (Note: here we need complete type-spaces.)
\end{lemma}

\begin{proof}
We may assume $B$ is $\l (\K)$-saturated.
For each $i<\o$, let $\d_i$ be as in Lemma~\ref{lemma5} (ii). By Lemmas \ref{lemma12} and \ref{lemma11}
we can find an increasing sequence of finite sets $A_i$ such that
$a\da^{\d_i}_{A_i}A$ and  $\bigcup_{i<\o}A_i=A$. Further choose
$a_i$ such that $Lstp(a_i/A_i)=Lstp(a/A_i)$ and $a_i\da_{A_i}B$.
By Lemma~\ref{lemma15}, $a_{i+1}\da^{\d_i}_{A_i}B$ and thus
$d^{p}_{a}(Lstp(a_{i+1}/B),Lstp(a_i/B))\leq\d_i$ which implies
$d^{p}(t^g(a_{i+1}/B),t^g(a_i/B))\leq\d_i$. Now the types $t^g(a_i/B)$
form a $d^{p}$-Cauchy
sequence and thus have a limit $a'$, i.e., for any $\e>0$ there is $i<\o$ such that $d^{p}(t^g(a_i/B),t^g(a'/B))<\e$ and, as $B$ was $\l (\K)$-saturated, $d^{p}_{a}(Lstp^w(a_i/B),Lstp^w(a'/B))<\e$. 
If $A'\subset A$ is finite and $\e>0$ there is $n<\omega$ such that $A'\subseteq A_n$ and $d^p_a(Lstp^w(a'/B),Lstp^w(a_n/B))<\e$ and as $a_n\models Lstp(a/A_n)$, 
$d^p_a(Lstp(a'/A'),Lstp(a/A'))<\e$. As $\e>0$ was arbitrary we must have $d^p(Lstp(a'/A')=Lstp(a/A'))$ and this must hold for all finite $A'\subset A$. So $Lstp^w(a'/A)=Lstp^w(a/A)$. Finally  $a'\da^{0}_AB$ is witnessed by the pairs $(a_n,A_n)$. 
\end{proof}

\begin{corollary}\label{corollary20}
If $Lstp^{w}(a/A)=Lstp^{w}(b/A)$,
then $Lstp(a/A)=Lstp(b/A)$.
\end{corollary}

\begin{proof}
 For finite $A$ the claim is trivial,
so let $A$ be infinite. By Corollary~\ref{corollary13} there is a countable
$B\subseteq A$ such that $ab\da^{0}_BA$. Let $\A\supset A$ be
$\l(\K)$-saturated. By Lemma~\ref{lemma19} choose $a'$ and $b'$ such that
$Lstp^w(a'b'/B)=Lstp^w(ab/B)$ and $a'b'\da^{0}_B\A$. By Lemma~\ref{lemma18},
$Lstp^w(a'b'/A)=Lstp^w(ab/A)$ and thus, by homogeneity, they
have the same Galois-type. So there exists an automorphism $F$
such that $F\raj A=id$ and $F(a'b')=ab$. Then $ab\da^{0}_BF(\A)$ so
by Lemma~\ref{lemma18}, $Lstp^w(a/F(\A))=
Lstp^w(b/F(\A))$. As $F(\A)\supseteq A$ and $F(\A)$ is $\l(\K)$-saturated, $Lstp(a/A)=Lstp(b/A)$. 
\end{proof}

\begin{lemma}\label{lemma21}
Let $A$ be finite or countable. Then $a\da^{0}_AB$ if and only if $a\da_AB$.
\end{lemma}

\begin{proof}
 First assume $A$ is finite. Then the direction from right to left is Lemma~\ref{lemma11}(vi). To prove the claim from left to right, let $a\da^{0}_AB$ and $b$ be such that $Lstp(b/A)=Lstp(a/A)$ and $b\da_AB$. By Lemma~\ref{lemma11}(vi) $b\da^{0}_AB$ so by Lemma~\ref{lemma18} and Corollary~\ref{corollary20}, $Lstp(a/B)=Lstp(b/B)$ and thus $a\da_AB$.

Then assume $A$ is countable and $a\da^{0}_AB$. By Lemmas \ref{lemma19} and \ref{lemma18} we may assume $B$ is $\l(\K)$-saturated and $B\supset A$. Now if $a\nda_AB$, $t^g(a/B)$ splits strongly over $A$. 
So there are $b,c\in B$, some $\e>0$ and finite $A'\subset A$
satisfying $Lstp(b/A)=Lstp(c/A)$ but $d^p(t^g(b/A'a),t^g(c/A'a))>\e$.
We claim that then $a\nda^{2(\e)}_AB$, namely for any finite
$A''$ with $A'\subseteq A''\subset A$ if $Lstp(a'/A'')=Lstp(a/A'')$
and $a'\da_{A''}B$ we must have $d^p_a(Lstp^w(a'/B),Lstp^w(a/B))>2(\e)$.
Otherwise we would have $d^p(t^g(bA'a/\empty),t^g(bA'a'/\empty))\leq2(\e)$
and $d^p(t^g(cA'a/\empty),t^g(cA'a'))\leq2(\e)$ and by
stationarity $t^g(bA'a'/\empty)=t^g(cA'a'/\empty)$,
adding up to $d^p(t^g(bA'a/\empty),t^g(cA'a/\empty))\leq\e$, a contradiction.

For the other direction assume $a\da_AB$. By Lemma~\ref{lemma19} and Corollary~\ref{corollary20} let $Lstp(a'/A)=Lstp(a/A)$, $a'\da^{0}_AB$. By the previous direction $a'\da_AB$ and thus by stationarity, $t^g(a'/B)=t^g(a/B)$, i.e. $a\da^{0}_AB$. 
\end{proof}

Now combining Corollary~\ref{corollary13} and Lemma~\ref{lemma21} we get:

\begin{corollary}\label{corollary22}
For all $A$ and $a$, there is countable $B\subseteq A$
such that $a\da_{B}A$.
\end{corollary}

This should be compared to what one gets from mere stability. In \cite{hysh} it is shown that 
for stable homogeneous AEC's
there is
$\k (\K )<\beth_{(2^{LS(\K )})^{+}}$ such that for
all $a$ and $\l (\K )$-saturated $\A$ there is
$A\subseteq\A$ of power $<\k (\K )$ such that
$a\da_{A}\A$. Even in the first-order case $\k (\K )$
cannot be chosen to be smaller than $LS(\K )^{+}$. 

We finally have the ingredients for simplicity, which will guarantee transitivity, symmetry and strong extension for $\da$ over any set.

\begin{corollary}\label{corollary23}
$\K$ is simple, i.e. $a\da_AA$ holds for any $a$ and $A$.
\end{corollary}

\begin{proof}
 Follows from Corollary~\ref{corollary22} and monotonicity of $\da$. 
\end{proof}

Note that weak simplicity does not in general imply simplicity. An example of a class that that is homogeneous, stable and weakly simple but not simple can be constructed by modifying an example by Shelah in \cite{hylerank} showing that $\o$-stability does not imply simplicity in the setting of homogeneous models.

\begin{example}
Let the vocabulary contain a binary relation symbol $E_i$ for each $i<\o+\o$. We let our monster model $\M$ consist of functions $f:\o+\o\to\kappa$ such that for some $i<\o+\o$ for all $j>i$ $f(j)=0$. On this model we let $E_i$ be an equivalence relation such that  $(f,g)\in E_i$ if
\begin{itemize}
\item[(a)] $i<\o$ and $f\raj i+1=g\raj i+1$ or
\item[(b)] $i\geq\o$, $f\raj\o=g\raj\o$ and for all $j>i$ $f(j)=g(j)$.
\end{itemize}
Then the class consisting of elementary submodels of $\M$ is homogeneous
and stable. It is not simple: let, for $n<\omega$, $f_n$ be such
that $f_n(i)=1$ if $i\leq n$ and $f_n(i)=0$ otherwise and define
$A=\{f_n:n<\omega\}$. Further let $f$ be such that $f(i)=1$ if $i<\o$
and $f(i)=0$ otherwise. Then $t^g(f/A)$ has no free extension so $f\nda_AA$.
However, this is the only way we do not get free extensions, so the
class is weakly simple.
\end{example}

\begin{lemma}\label{lemma24}
For every $\e>0$ there exists some $\d>0$ such that if $d^p_a(Lstp^w(a/A),Lstp^w(b/A))\leq\d$ and $ab\da^{0}_AB$, then $d^p_a(Lstp^w(a/B),Lstp^w(b/B))\leq\e$.
\end{lemma}

\begin{proof}
 Let $\d=3(\e)$. First note that by transitivity of
$\da^{0}$ we may assume $A$ to be countable. Let $B'\subset B$
be finite. We need to show $d^p_a(Lstp(a/B'),Lstp(b/B'))\leq\e$.
For this let $D\supset B'$ be finite and such that $D\da_{B'}ab$.
We may assume $D\da_{B'}abA$ and thus $ab\da_{B'A}D$. Now by Lemma~\ref{lemma21}
$ab\da^{0}_{B'A}D$ and thus by Lemma~\ref{lemma17} $ab\da^{0}_AB'D$. Now let $A'\subset A$
be finite and such that $ab\da^{\d}_{A'}D$ and choose $a'b'$ satisfying
$Lstp(a'b'/A')=Lstp(ab/A')$ and $a'b'\da_{A'}D$.
Then $d^p_a(Lstp(a/D),Lstp(a'/D))\leq\d$ and $d^p_a(Lstp(b'/D),Lstp(b/D))\leq\d$ and by Lemma~\ref{lemma8} (iii), $d^p_a(Lstp(a'/D),Lstp(b'/D))\leq\d$. This sums up to $d^p(t^g(a/D),t^g(b/D))\leq\e$. 
\end{proof}

\begin{lemma}\label{lemma25}
For every $\e>0$ there is a $\d>0$ such that if $a\da^{\d}_AB$ and $a\da_{AB}C$ then $a\da^{\e}_AC$.
\end{lemma}

\begin{proof}
 Let $\d$ be given by Lemma~\ref{lemma24} and assume $a\da^{\d}_AB$ and
$a\da_{AB}C$. Note that we may assume $A\subseteq B\subseteq C$ and
by Corollary~\ref{corollary22} and transitivity of $\da$ we may assume $B$ is countable.
Now if $a\nda^{\e}_AC$, then there is some finite $A'\subseteq A$ such
that $a\nda^{\e}_{A^+}C$ for all finite $A^+$ with $A'\subseteq A^+\subseteq A$.
As $a\da^{\d}_AB$, by Lemma~\ref{lemma11}, there is some finite $A''$ with
$A'\subseteq A''\subseteq A$ such that $a\da^{\d}_{A''}B$ and we still have
$a\nda^{\e}_{A''}C$. So
if there is a counterexample to the claim, we may find one with $A$ finite, $B$ at most countable and $A\subseteq B\subseteq C$ so it is enough to prove the lemma for such sets.

To prove $a\da^{\e}_AC$, let $b$ be such that $Lstp(b/A)=Lstp(a/A)$ and $b\da_AC$. We may assume that $b\da_ACa$. As $a\da^{\d}_AB$, we have $d^p_a(Lstp^w(a/B),Lstp^w(b/B))\leq\d$. Also now $ab\da_BC$ so by Lemma~\ref{lemma21} $ab\da^{0}_BC$ and thus by Lemma~\ref{lemma24} $d^p_a(Lstp^w(a/C),Lstp^w(b/C))\leq\e$. 
\end{proof}

\begin{corollary}\label{corollary26}
For every $\e>0$ there is a $\d>0$ such that if $a\da^{\d}_AB$ then for all $C$ there is $b$ satisfying $Lstp(b/AB)=Lstp(a/AB)$ and $b\da^{\e}_ABC$.
\end{corollary}

\begin{proof}
 Follows from Lemma~\ref{lemma25}
 by taking $b\models Lstp(a/AB)$ satisfying $b\da_{AB}C$. 
\end{proof}

We have seen that $d^p$-superstability and weak simplicity imply simplicity. We have also seen that $\da$ and $\da^{0}$ agree over finite and countable sets, and that $\da^{0}$ satisfies many of the properties of a well-behaved independence notion (left and right monotonocity, local character, transitivity, extension and stationarity of Lascar types). We still lack base monotonicity and symmetry and to prove them we look at a characterisation of $\da^{\e}$ using Lascar splitting.

\section{Lascar $\e$-splitting}\label{sec:esplit}

In this section we define and study Lascar $\e$-splitting. Via a characterisation of $\da^{\e}$ using Lascar splitting we can finally prove monotonicity of $\da^{0}$ and show that $\da$ and $\da^{0}$ are equal over all sets.

\begin{definition}\label{definition27}
\begin{itemize}
\item[(i)] If $A$ is finite and $A\subseteq B$, we say that $t^g(a/B)$ Lascar $\e$-splits over $A$ if for all $\d >0$, there are tuples $b,c\in B$
such that $d^{p}_{a}(Lstp(b/A),Lstp(c/A))<\d$ but
$d^{p}(t^g(ab/A),t^g(ac/A))>\e$.
\item[(ii)] We say that $t^g(a/B)$ locally Lascar $\e$-splits over $A\subseteq B$ if it Lascar $\e$-splits over every finite $A'\subseteq A$. (Note that for finite $A$ this is equivalent to (i).)
\end{itemize}
\end{definition}

\begin{lemma}\label{lemma28}
For all $\e >0$, there is $\d >0$ such that
if $a\da^{\d}_{A}B$, $A\subseteq B$, then
$t^g(a/B)$ does not locally Lascar $\e$-split over $A$.
\end{lemma}

\begin{proof}
 Let $\d =3(\e )$, let $A'\subseteq A$ be finite such that $a\da^{\d}_{A'}B$, and let $a'$ be such that $Lstp(a'/A')=Lstp(a/A')$
and $a'\da_{A'}B$. Now let $b,c\in B$ be such that
$d^{p}_{a}(Lstp(b/A'),Lstp(c/A'))<\d$. It suffices to show that
$d^{p}(t^g(ab/A'),t^g(ac/A'))\le\e$.
Since $a\da_{A'}^{\d}B$,
$d^{p}(t^g(ab/A'),t^g(a'b/A'))$, $d^{p}(t^g(ac/A'),t^g(a'c/A'))<\d$.
By the choice of $b$ and $c$, and by Lemma~\ref{lemma8} (iii), $d^{p}(t^g(a'b/A'),t^g(a'c/A'))<\d$. By Lemma~\ref{lemma9} this gives the required distance.
\end{proof}

\begin{thm}\label{theorem29}
For $A\subseteq B$ the following are equivalent:
\begin{itemize}
\item[(i)] $a\da_{A}B$,
\item[(ii)] for all $\e >0$, there is a finite $C\subseteq A$ with the following property:
for all $B'\supseteq B$ there is
$b$ such that $t^g(b/B)=t^g(a/B)$ and $t^g(b/B')$ does not
Lascar $\e$-split over $C$.
\end{itemize}
\end{thm}

\begin{proof}
 (i)$\Rightarrow$(ii): Let $\e>0$ be given. Let $\d$ be as in Lemma~\ref{lemma28} for $\e$ and let $\d'$ be as in Lemma~\ref{lemma25} for $\d$. Then choose a finite $C\subseteq A$ so that $a\da^{\d'}_{C}A$. Now for any $B'\supseteq B$ there is $b$ such that $t^g(b/B)=t^g(a/B)$ and $b\da_BB'$. Then $b\da_AB'$ and by Lemma~\ref{lemma25}, $b\da^{\d}_CB'$. By Lemma~\ref{lemma28} we are done.

(ii)$\Rightarrow$(i): Let $\D\supseteq B$ be a saturated model
of power $>\vert B\vert$. 
For all $n>0$, choose $b_{n}$ and $C_{n}$ as in (ii) for
$\e =1/n$ and $B'=\D$. We can choose these so that in addition,
for all $n>0$, $Lstp(b_{n}/B)=Lstp(a/B)$:
By the choice of $b_{n}$ there is an automorphism $F$ such that
$F\raj B=id$ and $F(b_{n})=a$. Then choose $b'\in F(\D )$
and $b''\in\D$ so that
$Lstp(b''/B)=Lstp(b'/B)=Lstp(a/B)$. Since $\D$ and $F(\D )$ are saturated,
there is an automorphism $G$ such that $G(F(\D ))=\D$,
$G\raj B=id$ and $G(b')=b''$. Now $G(a)$ is as wanted.

Let $c$ be such that $Lstp(c/A)=Lstp(a/A)$ and
$c\da_{A}\D$ and let $d\in B$. It is enough to show that
$t^g(cd/\empty )=t^g(ad/\empty )$. Let $\e >0$. It is enough
to show that $d^{p}(t^g(cd/\empty ),t^g(ad/\empty ))<\e$.
Let $n>0$ be such that $1/n <\e$. Since $d\in B$,
it is enough to show that $d^{p}(t^g(cd/\empty ),t^g(b_{n}d/\empty ))<\e$.
Choose $d'\in\D$ so that $Lstp(d'/A)=Lstp(d/A)$
and $d'\da_{A}b_{n}c$. Then
$t^g(cd/\empty)=t^g(cd'/\empty)=t^g(b_{n}d'/\empty )$.
Thus it is enough to show that
$d^{p}(t^g(b_{n}d/\empty ),t^g(b_{n}d'/\empty ))<\e$.
But since $t^g(b_{n}/\D)$ does not Lascar $\e$-split over
$C_{n}\subseteq A$, even
$d^{p}(t^g(b_{n}d/C_{n}),t^g(b_{n}d'/C_{n}))\le 1/n<\e$. 
\end{proof}

\begin{definition}\label{definition30}
We say that $A$ is almost strongly
$\o$-saturated if for all finite $B\subseteq A$, $\e >0$ and
$a$ there is $b\in A$ such that $d^{p}_{a}(Lstp(b/B),Lstp(a/B))<\e$.
\end{definition}

\begin{corollary}\label{corollary31}
Galois-types over almost strongly $\o$-saturated sets are
stationary.
\end{corollary}

\begin{proof}
 Let $A$ be almost strongly $\o$-saturated, $A\subseteq B$,
$t^g(a/A)=t^g(b/A)$, $a\da_AB$ and $b\da_AB$. Now if
$t^g(a/B)\neq t^g(b/B)$ there is $\e>0$ and some finite $B'\subset B$
such that $d^p(t^g(a/B'),t^g(b/B'))>\e$. Let $\d<2(\e)$.
By Theorem~\ref{theorem29} there is some finite $A'\subset A$ such that
$t^g(a/B)$ and $t^g(b/B)$ do not Lascar $\d$-split over $A'$, so
for some $\d'>0$ whenever $c_1,c_2\in B$ satisfy
$d^p_a(Lstp(c_1/A'),Lstp(c_2/A'))<\d'$ we have
$d^p(t^g(ac_1/A'),t^g(ac_2/A'))\leq\d$ and
$d^p(t^g(bc_1/A'),t^g(bc_2/A'))\leq\d$.
Now let $B^{\d'}\subset A$ satisfy $d^p_a(Lstp(B^{\d'}/A'),Lstp(B'/A'))<\d'$. Then $d^p(t^g(aB^{\d'}/A'),t^g(aB'/A'))\leq\d$ and $d^p(t^g(bB^{\d'}/A'),t^g(bB'/A'))\leq\d$. As $t^g(a/B^{\d'})=t^g(b/B^{\d'})$ this gives $d^p(t^g(a/B'),t^g(b/B'))\leq\e$, a contradiction. 
\end{proof}

By taking a closer look at the proof of (ii)$\Rightarrow$(i)
from Theorem~\ref{theorem29}, we get the following:

\begin{thm}\label{theorem32}
For all $\e >0$ there is $\d >0$ such that
if $B\supseteq A$ then for all $a$,
(*) below implies that $a\da^{\e}_{A}B$.
\begin{itemize}
\item[(*)] For all $D\supseteq B$ there is $b$ such that
$t^g(b/B)=t^g(a/B)$ and $t^g(b/D)$ does not locally Lascar $\d$-split over
$A$.
\end{itemize}
\end{thm}

\begin{proof}
 We first prove that (*) implies the following:

(*') There is a finite $A'\subset A$ such that for all $D\supseteq B$ there is $b$ such that $t^g(b/B)=t^g(a/B)$ and $t^g(b/D)$ does not Lascar $\d$-split over $A'$.

For this assume (*) and let $\D\supset B$ be a saturated model of power $>\abs{B}$. 
By (*) there is $b$ such that $t^g(b/B)=t^g(a/B)$ and $t^g(b/\D)$ does
not Lascar $\d$-split over some finite $A'\subset A$.
Then let $D'\supset B$ be any set and choose $b'$ such
that $t^g(b'/B)=t^g(a/B)$ and $b'\da_BD'$. We claim that $t^g(b'/D')$
does not Lascar $\d$-split over $A'$. Otherwise for any $\d'>0$
there are $c,d\in D'$ such that $d^p_a(Lstp(c/A'),Lstp(d/A'))<\d'$
but $d^p(t^g(b'c/A'),t^g(b'd/A'))>\d$. Then there is an automorphism
$F$ such that $F\raj B=id$ and $F(b')=b$. Denote $c'=F(c)$, $d'=F(d)$.
As $b'\da_BD'$, we have $b\da_Bc'd'$ and we may assume $c'd'\da_Bb\D$. By saturation of $\D$ we can find $c^+,d^+\in\D$ such that $Lstp(c^+d^+/B)=Lstp(c'd'/B)$ and $c^+d^+\da_Bb$. Thus $t^g(c^+d^+/Bb)=t^g(c'd'/Bb)$ and $t^g(c^+d^+b/B)=t^g(c'd'b/B)=t^g(cdb'/B)$. In particular $d^p_a(Lstp(c^+/A'),Lstp(d^+/A'))<\d'$ and $d^p(t^g(bc^+/A'),t^g(bd^+/A'))>\d$, a contradiction.

Now for the theorem, let $\d =2(\e )$. To prove $a\da^{\e}_AB$, let $C\subseteq A$ be finite. Let $A'\subseteq A$ be finite as given by (*') and define $A^+=A'\cup C$. Then let $b$ satisfy $Lstp(b/A^+)=Lstp(a/A^+)$ and $b\da_{A^+}B$.
Let $\D\supseteq B$ be a large saturated model
such that $\D\da_{B}ab$ (in particular $b\da_{A^+}\D$).
Clearly it is enough to show that $d^{p}(t^g(b/\D ),t^g(a/\D ))\le\e$.
For this, let $d\in\D$. It is enough to show that
$d^{p}(t^g(bd/\empty ),t^g(ad/\empty ))\le\e$. Let $b'$
be such that $t^g(b'/B)=t^g(a/B)$ and
$t^g(b'/\D )$ does not Lascar $\d$-split over $A'$ (and thus not over $A^+$).
As in the proof of (ii)$\Rightarrow$(i)
in Theorem~\ref{theorem29}, we can choose $b'$ so that in addition
$Lstp(b'/B)=Lstp(a/B)$.
By Lemma~\ref{lemma5}, it is enough to show that
$d^{p}(t^g(bd/\empty ),t^g(b'd/\empty ))\le\d$ and
$d^{p}(t^g(b'd/\empty ),t^g(ad/\empty ))\le\d$. 

For the first one choose $d'\in\D$ such that $Lstp(d'/A^+)=Lstp(d/A^+)$ and $d'\da_{A^+}bb'$. Then $t^g(bd/\empty)=t^g(bd'/\empty)=t^g(b'd'/\empty)$ and as $t^g(b'/\D)$ does not Lascar $\d$-split over $A^+$, $d^p(t^g(b'd/\empty),t^g(b'd'/\empty))\leq\d$.

For the second, choose $d'\in\D$ so that $Lstp(d'/B)=Lstp(d/B)$ and $d'\da_{B}ab'$.
Then $t^g(ad/\empty )=t^g(ad'/\empty )=t^g(b'd'/\empty )$
and since $t^g(b'/\D )$ does not Lascar $\d$-split over $A^+\subset B$,
$d^{p}(t^g(b'd'/\empty ),t^g(b'd/\empty ))\le\d$. 
\end{proof}

\begin{corollary}\label{corollary33}
For all $\e >0$ there is $\d >0$
such that for all $A\subseteq B\subseteq C$ and $a$,
if $a\da^{\d}_{A}C$,
then $a\da^{\e}_{B}C$.
\end{corollary}

\begin{proof}
 Let $\d=3(3(2(\e)))$ and assume $a\da^{\d}_AC$. To use Theorem~\ref{theorem32}, let $D\supseteq C$ and choose $b$ such that $t^g(b/C)=t^g(a/C)$ and $b\da_CD$. Then by Lemma~\ref{lemma25}, $b\da^{3(2(\e))}_AD$. By Lemma~\ref{lemma28}, $t^g(b/D)$ does not locally Lascar $2(\e)$-split over $A$ and thus not over $B$. By Theorem~\ref{theorem32}, $a\da^{\e}_BC$. 
\end{proof}

\begin{corollary}\label{corollary34}
For all $\e >0$, there is $\d >0$ for which
there are no $a$ and $b_{n},c_{n}$, $n>0$,
such that for all $n>0$, the following holds:
\begin{itemize}
\item[(i)] $d^{p}_{a}(Lstp(b_{n}/A_{n}),Lstp(c_{n}/A_{n}))<\d$,
where $A_{n}=\bigcup_{i<n}b_{i}c_{i}$,
\item[(ii)] $d^{p}(t^g(ab_{n}/A_{n}),t^g(ac_{n}/A_{n})) >\e$.
\end{itemize}
\end{corollary}

\begin{proof}
 Let $\d=3(3(2(3(\e ))))$.
For a contradiction, suppose that
$a$, $b_{n},c_{n}$ for $n<\o$ exist
such that (i) and (ii) hold.
We can find a finite $A\subseteq\bigcup_{n<\o}A_{n}$
such that $a\da^{\d}_{A}\bigcup_{n<\o}A_{n}$.
Choose $n<\o$ so that $A\subseteq A_{n}$. By the proof of Corollary~\ref{corollary33},
$a\da^{3(\e )}_{A_{n}}b_{n}c_{n}$. Let $b$ satisfy $Lstp(b/A_n)=Lstp(a/A_n)$ and $b\da_{A_n}b_nc_n$. Then $d^p_a(Lstp(b/A_nb_nc_n),Lstp(a/A_nb_nc_n))\leq 3(\e)$.
Further since $d^{p}_{a}(Lstp(b_{n}/A_{n}),Lstp(c_{n}/A_{n}))<\d\le 3(\e)$, by Lemma~\ref{lemma8} (iii) we also have that $d^p_a(Lstp(b_n/A_nb),Lstp(c_n/A_nb))\leq 3(\e)$. Then finally by Lemma~\ref{lemma9} this sums up to $d^p_a(Lstp(ab_n/A_n),Lstp(ac_n/A_n))\leq\e$, a contradiction. 
\end{proof}

Corollary~\ref{corollary33} ensures full monotonicity for $\da^{0}$, giving us the following generalisations of Lemmas \ref{lemma21} and \ref{lemma19}:

\begin{corollary}\label{corollary35}
$a\da^{0}_AB$ if and only if $a\da_AB$.
\end{corollary}

\begin{proof}
 Assume $a\da^{0}_AB$. By Corollary~\ref{corollary13} let $A_0\subseteq A$ be countable and such that $a\da^{0}_{A_0}A$. By Lemma~\ref{lemma17}, $a\da^{0}_{A_0}B$ and by Lemma~\ref{lemma21}, $a\da_{A_0}B$, and thus $a\da_AB$.

For the other direction suppose $a\da_AB$. Again let $A_0\subseteq A$ be countable and such that $a\da^{0}_{A_0}A$. By Lemma~\ref{lemma25}, $a\da^{0}_{A_0}B$ and by Corollary~\ref{corollary33}, $a\da^{0}_AB$. 
\end{proof}

\begin{corollary}\label{corollary36}
For any $a$ and $A\subseteq B$ there exists some $a'$ satisfying  $Lstp(a'/A)=Lstp(a/A)$ and $a'\da^{0}_{A}B$. 
\end{corollary}

\section{Almost summability and $\da^{>\e}$}\label{sec:almost}

In this section we study a property which we call almost summability. It allows us to add small distances to a given distance in the type space without the combined distance growing too big. We also study a weakening of $\da^{\e}$, which under the assumption of almost summability is very well behaved.

\begin{definition}\label{definition37}
We say that the perturbation system $(\F_\e)_{\e\geq0}$ is \emph{almost summable} if for all $\e>\d>0$ there exists some $m(\e,\d)>0$ such that for all $a_i$, $i\leq 2$, if $d^p(t^g(a_0/\emptyset),t^g(a_1/\emptyset))\leq\d$ and $d^p(t^g(a_1/\emptyset),t^g(a_2/\emptyset))\leq m(\e,\d)$ then $d^p(t^g(a_0/\emptyset),t^g(a_2/\emptyset))\leq\e$.
\end{definition}

\begin{remark}\label{remark38}
\begin{itemize}
\item[(i)] As in Lemma~\ref{lemma9} one can show that if the perturbation system is almost summable then with $\e>\d>0$ and $m(\e,\d)$ as in the definition, also $d^p_a$-distances of $\d$ and $m(\e,\d)$ add up to $\e$.
\item[(ii)] Almost summability holds e.g.\ for the perturbation system of Hilbert spaces with an automorphism \cite{byusza} or linear isomorphisms of Banach spaces. Below we give an example where almost summability fails.
\end{itemize}
\end{remark}

\begin{example}\label{example39}
We give an example of a class that is homogeneous with complete type spaces 
but whose perturbation system is not almost summable. The vocabulary is
$L=\{P_n,E,<,R_q,d\}_{n<\omega,q\in\Q\cap(0,2]}$ where the $P_n$
are unary predicates and $E$, $<$ and $R_q$ are binary. $E$ is an
equivalence relation, the predicates $P_n$ partition the universe and
each predicate is a union of $E$-equivalence classes. $<$ is an order
on each equivalence class such that there for each equivalence class
exists a real $1\leq r\leq10$
such that $([a]_E,<)$ is isomorphic to the ordered real interval $[r,2r]$.
$R_q(a,b)$ holds if and only if $[a]_E=[b]_E$ and $b/a=q$.
The metric $d$ is defined as the one induced by the interval
$[r,2r]$ within the equivalence classes and $d(a,b)=10$ if
$a$ and $b$ are in different equivalence classes. $d$ and
$<$ together fix the $r$ and a unique isomorphism $l:[a]_E\to[r,2r]$
for each element $a$. Thus we can define $r_a$ as the real $r$ given by
the isomorphism above and the length of $a$ as $l(a)\in[r_a,2r_a]$.

We define the perturbation system as follows: $f\in\F_\e$ if $f$ is a $L\backslash\{d\}$-isomorphism and if $a\in P_n$ then also
$$
e^{-n\e}\leq\frac{l(f(a))}{l(a)}\leq e^{n\e}.
$$
The above condition makes sure that $\F_0=\bigcap_{\e>0}\F_\e$. As the $R_q$
prevents $\e$-isomorphisms from stretching the interval $[r,2r]$ the error
in metric arises from mapping equivalence classes onto each other and
thus switching the $r$. As $r$ varies between 1 and 10 this can only
increase distances to the 10-fold and thus $\e$-isomorphisms are bi-Lipschitz
with Lipschitz constant 10 (regardless of $\e$) so they are uniformly
continuous. The rest of the conditions of a perturbation system are trivial.

It is not hard to see that this gives a MAEC with perturbations that is homogeneous with JEP, AP, the perturbation property and complete type spaces.
The perturbation system, however, is not almost summable. If $\e>\d>0$ is such that $\e<4\d$ and $\d<2$
we show that no $\d'>0$ can suffice as the $m(\e,\d)$ in Definition~\ref{definition37}. So
let $\d'>0$ be given and choose $n$ such that $e^{n\d'}>10$.
Within $P_n$ let $a,b,c$ be elements in an equivalence class corresponding
to an interval $[r,2r]$ with $r<2$ and such that $a<b<c$ and $d(b,c)=\d$.
Then $d^p(t^g(ab/\emptyset),t^g(ac/\emptyset))=d(ab,ac)=\d$. Now we can
map $a,b,c$ with a $\d'$-isomorphism to elements $a',b',c'$ in an
equivalence class, inside $P_n$, corresponding to an interval $[r',2r']$
with $r'>9$. This shows that
$d^p(t^g(ac/\emptyset),t^g(a'c'/\emptyset))\leq\d'$. But as this is the only way we can map $a,b,c$ to that interval, we must have $d^p(t^g(ab/\emptyset),t^g(a'c'/\emptyset))\geq d(a'b',a'c')>4\d>\e$.
\end{example}

\begin{definition}\label{definition40}
We define $a\da^{>\e}_AB$ if $a\da^{\xi}_AB$ for all $\xi>\e$. Note that with this notation $a\da^{0}_AB$ if and only if $a\da^{>0}_AB$.
\end{definition}

\begin{remark}\label{remark41}
If $A$ is finite $a\da^{>\e}_AB$ if and only if $a\da^{\e}_AB$. This is easily seen via the observation that if $A$ is finite then $a\da^{\e}_AB$ says that $a$ is $\e$-$d^p_a$-close to the free extension of $Lstp(a/A)$ over $B$. With an almost summable perturbation system a similar characterisation holds for any $A$:
\end{remark}

\begin{lemma}\label{lemma42}
Assume the perturbation system is almost summable and $A\subseteq B$.
\begin{itemize}
\item[(i)] If $a\da^{>\e}_AB$, $Lstp^w(b/A)= Lstp^w(a/A)$ and $b\da_AB$, then we have $d^p_a(Lstp^w(a/B),Lstp^w(b/B))\leq\e$. 
\item[(ii)] If $Lstp^w(b/A)= Lstp^w(a/A)$, $b\da_AB$ and $d^p_a(Lstp^w(a/B),Lstp^w(b/B))\leq\e$ then $a\da^{>\e}_AB$.
\end{itemize}
\end{lemma}

\begin{proof}
 (i) Assume $a\da^{>\e}_AB$, $b\models Lstp^w(a/A)$ and $b\da_AB$. We prove that $d^p_a(Lstp^w(a/B),Lstp^w(b/B))\leq\xi$ for every $\xi>\e$. So let $\xi>\e$ be given and let $\xi>\e''>\e'>\e$, $\d^+=\min\{m(\xi,\e''),m(\e'',\e')\}$
and $\d=3(\d^+)$. By Lemma~\ref{lemma12} find a finite $A_\d\subseteq A$ such that $ab\da^{\d}_{A_\d}A$. As in the proof of Lemma~\ref{lemma17} (since $\d\leq m(\e'',\e')$) we get $a\da^{\e''}_{A_\d}B$. Now let $c\models Lstp(a/A_\d)$ and $c\da_{A_\d}Bb$. Then $d^p_a(Lstp^w(a/B),Lstp^w(c/B))\leq\e''$. Also, as $c\da_{A_\d}A$, we have $d^p_a(Lstp^w(b/A),Lstp^w(c/A))\leq\d$. Further by Fact~\ref{fact2} and Corollary~\ref{corollary35}, $bc\da^{0}_AB$, so by Lemma~\ref{lemma24} $d^p_a(Lstp^w(b/B),Lstp^w(c/B))\leq\d^+$. But then $d^p_a(Lstp^w(a/B),Lstp^w(b/B))\leq\xi$.

(ii) Assume $b\models Lstp^w(a/A)$, $b\da_AB$ and
$d^p_a(Lstp^w(a/B),Lstp^w(b/B))\leq\e$ and let $\xi>\e$.
By Corollary~\ref{corollary35} $b\da^{0}_AB$ so by Lemma~\ref{lemma11} for every finite $C\subseteq A$ there is a finite $A'$ with $C\subseteq A'\subseteq A$ such that $b\da^{m(\xi,\e)}_{A'}B$. Now if $c\models Lstp(a/A')=Lstp(b/A')$ and $c\da_{A'}B$, we have $d^p_a(Lstp^w(b/B),Lstp^w(c/B))\leq m(\xi,\e)$.
Together with the assumption this yields $d^p_a(Lstp^w(a/B),Lstp^w(c/B))\leq\xi$ proving $a\da^{\xi}_AB$. 
\end{proof}

\begin{corollary}\label{corollary43}
If the perturbation system is almost summable, $a\da_AB$ and $a\da^{>\e}_{AB}C$ then $a\da^{>\e}_ABC$.
\end{corollary}

\begin{proof}
Assume $a\da_AB$ and $a\da^{>\e}_{AB}C$ and let $b\models Lstp^w(a/A)$, $b\da_ABC$. By stationarity $Lstp^w(b/AB)=Lstp^w(a/AB)$ so by $a\da^{>\e}_{AB}C$ and Lemma~\ref{lemma42} (i) $d^p_a(Lstp^w(a/ABC),Lstp^w(b/ABC))\leq\e$. But then by Lemma~\ref{lemma42} (ii) $a\da^{>\e}_ABC$. 
\end{proof}

\begin{lemma}\label{lemma44}
If the perturbation system is almost summable then $a\da^{>\e}_ABC$ and $a\da_AB$ imply $a\da^{>\e}_{AB}C$.
\end{lemma}

\begin{proof}
Let $\xi>\e$ and some finite $B'\subseteq AB$ be given. Let $\xi>\e'>\e$, $\d=m(\xi,\e')$ and $\d'=3(3(2(\d)))$ (from Corollary~\ref{corollary33}). Denote $A'=B'\cap A$. By $a\da_AB$ there is a finite $A''\supseteq A'$ such that $a\da^{\d'}_{A''}AB$. By $a\da^{>\e}_ABC$ there is a finite $A^+\supseteq A''$ such that $a\da^{\e'}_{A^+}ABC$. Define $B^+=A^+\cup B'$ and let $b\models Lstp(a/B^+)$ such that $b\da_{B^+}ABC$. We need to show that $d^p_a(Lstp^w(a/ABC),Lstp^w(b/ABC))\leq\xi$.

Now by Corollary~\ref{corollary33} and the choice of $\d'$ we have $a\da^{\d}_{A^+}AB$ and thus $b\da^{\d}_{A^+}B^+$. By Lemma~\ref{lemma15} we get $b\da^{\d}_{A^+}ABC$. Let $b'\models Lstp(a/A^+)=Lstp(b/A^+)$ such that $b'\da_{A^+}ABC$. Then $d^p_a(Lstp^w(b/ABC),Lstp^w(b'/ABC))\leq\d$. By $a\da^{\e'}_{A^+}ABC$ we also have $d^p(Lstp^w(a/ABC),Lstp^w(b'/ABC))\leq\e'$. By almost summability we are done. 
\end{proof}

We will use the notion $\da^{>\e}$ when studying pregeometries and a related rank in section~\ref{sec:pregeom}. However, before that we have a look at another example, related to the U-rank.

\section{Entropy}\label{sec:entropy}

Our first example on the use of measures of dependence is entropy.
In \cite{BH}, Berenstein and Henson showed a connection between
entropy and $\e$-dividing in the
(continuous first-order) context of probability algebras.
Below we do much the same in our context but in order to avoid
just repeating what was done in
\cite{BH}, we change the point of view a bit.

In 1948, C. Shannon suggested that entropy can be seen as
a measure of information or uncertainty. He studied the question on
how to measure the uncertainty if we know just probabilities
$p_{1},\ldots,p_{n}$  of possible events. He suggested that this measure should
satisfy three very reasonable requirements and then went on to
show that the only functions that satisfy the requirements
are of the form
$$H(p_{1},\ldots,p_{n})=-K\sum_{i=1}^{n}p_{i}\log(p_{i}),$$
where $K\in\R_{+}$, see \cite{Pe}. From now on we will use $K=1$
and the reader is free to choose the logarithm.
One of Shannon's requirements was that for all
$0\le\l\le 1$, $H(p_{1},\ldots,p_{n-1},\l p_{n},(1-\l )p_{n})=
H(p_{1},\ldots,p_{n})+p_{n}H(\l ,(1-\l ))$. In particular,
$H(p_{1},\ldots,p_{n-1},\l p_{n},(1-\l )p_{n})\ge
H(p_{1},\ldots,p_{n})$, which is a usefull fact in calculations.

When Shannon was asked what he had thought about when he had
confirmed his measure,
he had answered (see \cite{TM}): My greatest concern was what to call it.
I thought of calling it information, but
the word was overly used, so I decided to call it uncertainty.
When I discussed with John von Neumann, he had a better idea.
``You should call it entropy, for two reasons. In the first place
your uncertainty has been used in statistical mechanics
under that name, so it already has a name. In the second place, and more
importantly, no one knows what entropy really is,
so in a debate you will always have the advantage.''

Since types are packages of information, it may make sense
to define entropy for them. This can be done as follows:
Suppose $A\subseteq B$. By $E(a,A,B)$ we denote
the supremum of all $\e>0$ such that $a\nda^{\e}_{A}B$.
If there are no such $\e$, i.e., $a\da_{A}B$,
we let $E(a,A,B)=0$. Since we are going to need scaling,
for any  strictly increasing continuous  function
$\pi$ from $\R_{+}$ onto $\R_{+}$, by
$E_{\pi}(a,A,B)$ we denote $\pi (E(a,A,B))$.
The entropy $H_{\pi}(a/A)$ of the type $t(a/A)$ is then defined as
supremum of
$$-\sum_{i=1}^{n}E_{\pi}(a,A_{i-1},A_{i})\log(E_{\pi}(a,A_{i-1},A_{i}))$$
over all 
$0<n<\o$ 
and $A_{i}$, $i\le n$, such that
$A_{0}=A$ and $A_{i}\subseteq A_{i+1}$
and here $0\log(0)$ is defined to be $0$ (notice that
$\lim x\log(x)=0$ when $x$ goes to zero).

Let us look at how this notion behaves in a very basic setup
of quantum mechanics: We let $\K$ be the class of all
models of the form $(H,B)$, where $H$ is a Hilbert space over the
complex numbers and $B\subseteq H$ is an orthonormal basis
of $H$. We let
$\preceq$ be the submodel relation and we let
the perturbations be trivial, i.e.,
$\BF_{\e}=\BF_{0}$ for all $\e >0$.
Notice that this class is not axiomatizable in the continuous
first-order logic since it is not closed under ultraproducts
(the ultraproduct of $B$ is not a basis for the ultraproduct of $H$).
Notice also that if $u\in H$, $w\in B$ and
$\langle u,w\rangle\ne 0$, then $w\in bcl(\{ u\} )$.

Now let $(H,B)\in\K$ be such that $H$ is separable (not necessarily of
infinite dimension). We think of $B=\{ v_{i}\vert\  i<N\}$, $N\le\o$,
as a set of eigenvectors with eigenvalues $\l_{i}$ for some
observable $P$ such that $\l_{i}\ne\l_{j}$ for $i\ne j$.
For simplicity we leave $P$ out of the models, since 
the
exact
eigenvalues themselves
do not play a role in Shannon's entropy.
Now (e.g)
let $v=\sum_{i=0}^{n}a_{i}v_{i}$ be a state, where
$n\le N$ and $n<N$ if $N=\o$. Then the probability
for the observable getting the value $\l_{i}$ in a measurement is
$p_{i}=\abs{a_{i}}^{2}$ if $i\le n$ and is $0$ otherwise.
Now letting our scaling function $\pi$ be such that
$\pi (x)=x^{2}/2$, then we can calculate as follows
(easy, we leave the details to the reader):
If for all $i\le n+1$, we let $A_{i}=\{ v_{j}\vert\ j<i\}$, then
$$H_{\pi}(v/\empty)
=\sum_{i=1}^{n+1}E_{\pi}(a,A_{i-1},A_{i})\log(E_{\pi}(a,A_{i-1},A_{i}))=$$
$$\sum_{i=1}^{n+1}\abs{a_{i-1}}^{2}\log(\abs{a_{i-1}}^{2})=
H(p_{0},\ldots,p_{n}).$$

Note that by $\distp$-superstability, we can define an `$\varepsilon$-U-rank' counting the length of $\varepsilon$-forking chains. The ordinary U-rank need not be defined, as this would require discrete superstability (we may have infinite forking chains with decreasing measures of dependence). Entropy, however, can be defined even with infinite forking chains as long as the measure of dependence decreases fast enough, so it gives another rank function.

\section{Finding a pregeometry in $\M^{eq}$}\label{sec:pregeom}

In this section we study a closure operator defined by $\nda$ on the set of realisations of a Lascar type. We find conditions on $\nda^{\e}$ that guarantee that there is an equivalence relation on this set such that the closure operator forms a pregeometry on the set of equivalence classes. The $p$-adic integers, studied at the end of this paper, form an example of a class where this happens, but where the type itself is not regular (which is what is usually needed to find pregeometries).

Let $D$ be the set of all realisations of some unbounded $p=Lstp^w(a/A)$. 
Let $E$ be an $A$-invariant equivalence relation.
Denote by $a^*$ the $E$-equivalence class of $a$. We define in $D/E$ a closure operator by 
$a^*\in cl(b_1^*,\dots,b_n^*)$ if for all $a'\in a^*$ and $b_i'\in  b_i^*$, $i=1,\dots,n$, $a'\nda_Ab_1'\dots b_n'$. For an arbitrary $B^*\subseteq D/E$ we define $a^*\in cl(B^*)$ if $a^*\in cl(B_0^*)$ for some finite $B_0^*\subseteq B^*$.

\begin{lemma}\label{lemma45}
$cl$ as defined above satisfies Steinitz' exchange property, i.e., if $a^*\in cl(b_1^*,\dots,b_n^*,c^*)\backslash cl(b_1^*,\dots,b_n^*)$ then $c^*\in cl(b_1^*,\dots,b_n^*,a^*)$.
\end{lemma}

\begin{proof}
 Assume $a^*\in cl(b_1^*,\dots,b_n^*,c^*)\backslash cl(b_1^*,\dots,b_n^*)$.
If there are $c'\in c^*,a'\in a^*, b_k'\in b_k^*$, $1\leq k\leq n$,
such that $c'\da_A b_1'\dots b_n'a'$ then we can form a Morley sequence
$c_i$, $i<\lambda(\K)$ such that $c_i\models Lstp(c'/A)$ and
$c_i\da_Aa'b_1'\dots b_n'\bigcup_{j<i}c_j$. Now for each $i<\lambda(\K)$
there is an automorphism $F_i\in Aut(\M/Ab_1',\dots b_n'a')$ mapping
$c_i$ to $c'$ and since it fixes the $b_k'$, $1\leq k\leq n$, and $a'$ it
must fix their equivalence classes setwise.

Now as $a^*\notin cl(b_1^*,\dots,b_n^*)$ there are $a''\in a^*$ and $b_k''\in b_k^*$, for $1\leq k\leq n$, such that $a''\da_Ab_1''\dots b_n''$. Now $F_i(a'')\in a^*$, $F_i(b_k'')\in b^*$ and $F_i(c_i)=c'\in c^*$ and as $a^*\in cl(b_1^*,\dots,b_n^*,c^*)$ we have $F_i(a'')\nda_AF_i(b_1'')\dots F_i(b_n'')c'$ so $a''\nda_Ab_1''\dots b_n'' c_i$ for every $i<\lambda(\K)$. As $a''\da_A b_1''\dots b_n''$ this implies $a''\nda_{Ab_1''\dots b_n''}c_i$. By symmetry $c_i\nda_{Ab_1''\dots b_n''}a''$ 
and further $c_i\nda_Ab_1''\dots b_n''a''$. But as the $c_i$ form a Morley sequence this implies $b_1''\dots b_n''a''\nda_{A\bigcup_{j<i}c_j}c_i$ for every $i<\lambda(\K)$ but this gives a strongly splitting chain of length $\lambda(\K)$, a contradiction. So we must have $c'\nda_A b_1'\dots b_n'a'$ for all $a'\in a^*$, $b_k'\in b_k^*$, $1\leq k\leq n$. 
\end{proof}

\begin{lemma}\label{lemma46}
Assume $A$ is finite or the perturbation system is almost summable.
If (*) below holds, then $(D/E,cl)$ is a pregeometry.
\begin{itemize}
\item[(*)] There is $\e>0$ such that for all $b\in D$ and
$B\subseteq D$ the following are equivalent: 
\begin{itemize}
\item[(i)] $b\nda_AB$
\item[(ii)] $b\nda^{>\e}_AB$
\item[(iii)] for all $c\in D$ there exists $b'\in b^*$ such that $c\da^{>\e}_{AB} b'$.
\end{itemize}
\end{itemize}
\end{lemma}

\begin{proof}
Monotonicity is clear, finite character was built into the definition
and exchange was proved in Lemma~\ref{lemma45} so all that remains is $cl(cl(B))=cl(B)$.
It is enough to consider the case where $c^*\in cl(b_1^*,\dots,b_n^*)$ and
$a^*\in cl(b_1^*,\dots,b_n^*,c^*)$ and show that $a^*\in cl(b_1^*,\dots,b_n^*)$.
So let $a'\in a^*,b_i'\in b_i^*$. We need to show $a'\nda_Ab_1'\dots b_n'$.
As all $c\in c^*$ satisfy $c\nda_Ab_1'\dots b_n'$ by (*) there exists
$c'\in c^*$ such that $a'\da^{>\e}_{Ab_1'\dots b_n'}c'$.
Now if $a'\da_Ab_1'\dots b_n'$ and $A$ is finite then by Lemma~\ref{lemma11} (vii) $a'\da^{>\e}_A b_1'\dots b_n'c'$ and by (*) $a'\da_A b_1'\dots b_n'c'$, a contradiction. 
If $d^p$ is almost summable then by Corollary~\ref{corollary43} 
we again get $a'\da_A b_1'\dots b_n'c'$. So $a'\nda_Ab_1'\dots b_n'$. 
\end{proof}

\begin{remark}\label{remark47}
Note that in the lemma above $p$ itself need not be regular as we will see in the $p$-adic example in section~\ref{sec:ex}. 
\end{remark}

\begin{lemma}\label{lemma48}
If $E$ is an equivalence relation such that (*) of Lemma~\ref{lemma46} holds, then $aEb$ implies $a\nda_Ab$.
\end{lemma}

\begin{proof}
 We first show that if (*) holds then $E$ has more than one equivalence class. Assume $aEb$ and let $B$ be such that $a\nda_AB$. Now by the equivalence of (i) and (iii) in (*) and since $a^*=b^*$ we get $b\nda_AB$. Thus if $c\da_AB$, $c$ is not in the $E$-class of $a$, so $E\neq D^2$.

Now if $aEb$, $a\da_Ab$ and $c\in D$ is arbitrary, choose $d\in D$ with $d\da_Abc$ then $t^g(d/Ab)=t^g(a/Ab)$ so $dEb$. Further $t^g(c/Ad)=t^g(b/Ad)$ so $cEb$ and $E$ has $D$ as its only equivalence class, a contradiction. 
\end{proof}

\begin{corollary}\label{corollary49}
If $E$ is an equivalence relation such that (*) of Lemma~\ref{lemma46} holds and either $A$ is finite or the perturbation system is almost summable, then $\nda_A$ is transitive on $D$ (and thus an equivalence relation) and (*) holds for this equivalence relation.
\end{corollary}

\begin{proof}
Assume towards a contradiction that $a\nda_Ab$ and $b\nda_Ac$ but $a\da_Ac$. As $b\nda_Ac$ there is $b'$ such that $b'Eb$ and $a\da^{>\e}_{Ac}b'$ and as $a\da_Ac$ we have $a\da^{>\e}_Acb'$ and thus $a\da_Ab'$. Further, as $b\nda_Aa$ there is $b''$ such that $b''Eb$ and $b'\da^{>\e}_{Aa}b''$. Now as $b'\da_Aa$ we get $b'\da^{>\e}_Aab''$ and thus $b'\da_Ab''$, contradicting Lemma~\ref{lemma48}.

Now $\nda_A$ forms an equivalence relation on $D$ such that each $\nda_A$-equivalence class is a union of $E$-equivalence classes. Thus (iii) with respect to $E$ implies (iii) with respect to $\nda_A$.

Now if for all $c\in D$ there is $b'$ such that $b'\nda_Ab$ and $c\da^{>\e}_{AB}b'$ then in particular this holds for $c=b$. Then if $b\da_AB$ we have $b\da^{>\e}_ABb'$ and by equivalence of (i) and (ii) $b\da_Ab'$, a contradiction. 
\end{proof}

The relation $\da^{\e}$ (or $\da^{\e}$) measures distances to free extensions. Another view is looking at how much a type can still fork.

\begin{definition}\label{definition50}
We define a real-valued rank function
$$
R(a/A)=\sup\{\e : a\nda^{\e}_AB \textrm{ for some }B\} 
$$
and $R(a/A)=0$ if the above set is empty.
\end{definition}

\begin{remark}\label{remark51}
\begin{itemize}
\item[(i)] $R(a/A)=0$ if and only if $t^g(a/A)$ is bounded.
\item[(ii)] $R$ is not in general monotone (i.e.\ it is not always the case that $R(a/AB)\leq R(a/A)$ as can be seen by considering functions $\omega\to X$ for some set $X$ and defining the following metric:
$$
d(f,g) = 
\left\{\begin{array}{ll}
0,9 &\textrm{if }\min\{n:f(n)\neq g(n)\}=0\\
m^{-1} &\textrm{if }m=\min\{n:f(n)\neq g(n)\}>0.\end{array}\right.
$$
\end{itemize}
\end{remark}

\begin{lemma}\label{lemma52}
If either $A$ is finite or the perturbation system is almost summable then $a\da_AC$ implies $R(a/AC)=R(a/A)$.
\end{lemma}

\begin{proof}
  Assume $a\da_AC$. If $R(a/A)>\e$ then for some $B$ $a\nda^{>\e}_AB$ and thus by right monotonicity $a\nda^{>\e}_ABC$. As $a\da_AC$ this implies $a\nda^{>\e}_{AC}B$ (by Lemma~\ref{lemma11}(vii) if $A$ is finite and by Corollary~\ref{corollary43} if the perturbation system is almost summable) which shows that $R(a/AC)>\e$.

If $R(a/A)\leq\e$ then for all $B$ $a\da^{\e}_AB$, in particular $a\da^{\e}_ABC$. Then either by Lemma~\ref{lemma11}(viii) or Lemma~\ref{lemma44} and $a\da_AC$ we get $a\da^{>\e}_{AC}B$ for all $B$, proving $R(a/AC)\leq\e$. 
\end{proof}

\begin{lemma}\label{lemma53}
Let $D$ be as above and assume either that $A$ is finite or the perturbation system is almost summable.  Suppose there is $\e>0$ such that for all $B\subset D$ and all $b\in D$ the following are equivalent
\begin{itemize}
\item[(i)] $b\nda_AB$
\item[(ii)] $b\nda^{>\e}_AB$
\item[(iii)] $R(b/AB)\leq\e$
\end{itemize}
Then $a\nda_Ab$ is an equivalence relation on $D$.
\end{lemma}

\begin{proof}
 Suppose $a,b,c\in D$, $b\nda_Aa$ and $b\nda_Ac$. For a contradiction suppose $a\da_Ac$. Then $b\nda_{Aa}c$. Let $b'\models Lstp^w(b/Aa)$ and $b'\da_{Aa}c$. Since $R(b/Aa)\leq\e$, we have $b\da^{>\e}_{Aa}c$ so $d^p_a(Lstp^w(b'/Aac),Lstp^w(b/Aac))\leq\e$. In particular, $d^p_a(Lstp^w(b'/Ac),Lstp^w(b/Ac))\leq\e$.
On the other hand, $b'\da_Ac$ and so by (ii), $d^p_a(Lstp^w(b'/Ac),Lstp^w(b/Ac))>\e$, a contradiction. 
\end{proof}

\begin{corollary}\label{corollary54}
Let $D$ be as above. Assume there is $\e>0$ such that for all $b\in D$ and $B\subseteq D$ the following are equivalent: 
\begin{itemize}
\item[(i)] $b\nda_AB$
\item[(ii)] $b\nda^{>\e}_AB$
\item[(iii)] $R(b/AB)\leq\e$
\item[(iv)] for all $c\in D$ there exists $b'\in D$ with $b'\nda_Ab$ such that $c\da^{>\e}_{AB} b'$
\end{itemize}
Then $a\nda_Ab$ is an equivalence relation on $D$ and $(D/\nda_A,cl)$ is a pregeometry.
\end{corollary}

\begin{proof}
This is clear by Lemmas \ref{lemma53} and \ref{lemma46}. 
\end{proof}

\section{Example: the $p$-adics}\label{sec:ex}

\def\ss{\preccurlyeq_{\K}}
\def\pcl#1{\left\langle#1\right\rangle_P}
\def\cl#1{\left\langle#1\right\rangle}
\def\pnorm#1{\left\Vert #1\right\Vert_p}
\def\div{\mid}
\def\ndiv{\nmid}
We finally demonstrate the properties studied in an example class of ultrametric spaces. The class consists of  models $\overline{\Z_p^{(\kappa)}}$ for a fixed prime $p$, where $\Z_p$ is the set of $p$-adic integers (i.e.\ the completion of the integers in the $p$-adic topology). Recall that the $p$-adic topology is given by the $p$-adic norm $\pnorm{a}=p^{-max\{k : p^k\div a\}}$.

We shortly recall some group theoretic notions. An element $a\in A$ is \emph{divisible by $n$} if there is $a'\in A$ such that $na'=a$. A subgroup $B\leq A$ is \emph{pure} if for every $n\in\Z$ each $b\in B$ which is divisible by $n$ in $A$ is divisible by $n$ already in $B$. The \emph{$p$-height} of $a$ is the largest $k\in\N$ such that $a$ is divisible by $p^k$. All $p$-adic integers are divisible by all $n\in\N$ coprime to $p$, so in the $p$-adic integers \emph{height} refers to $p$-height.

If $A$ is a subset of a group $B$, $\cl{A}$ denotes the subgroup generated by $A$ and $\pcl{A}$ denotes the pure subgroup in $B$ generated by $A$, i.e. $\pcl{A}=\{b\in B : \exists n\in\N^\ast\ nb\in\cl{A}\}$. When using this notation $B$ will be clear from the context.

We work in the vocabulary of groups $L=\{0,+,-\}$. The class $\K_p$ consists of completions of direct sums of copies of the $p$-adic integers, $\overline{\Z_p^{(\kappa)}}$ with $\kappa$ any cardinal. We let $A\ss B$ if $A$ is a closed pure subgroup of $B$.

Although we work in the vocabulary of groups, we will use the fact that our models, as completions in the $p$-adic topology of $\Z$-modules, are $p$-adic modules (modules over the ring of $p$-adic integers). Thus we can use structure theorems of complete modules. By a pure submodule of a $p$-adic module $A$ we mean a submodule $B$ such that $p^kB=B\cap p^kA$ for all $k\in\N$. Thus a pure closed subgroup of a $p$-adic module is a pure submodule.

The facts below actually hold for complete modules over any complete discrete valuation ring (complete principal ideal ring with exactly one prime element), but we constrain our attention to $p$-adic modules.

\begin{fact}[{\cite[\S 16]{kaplansky}}]\label{fact55}
\begin{itemize}
\item[(i)] The completion, in the $p$-adic topology, of a module with no elements of infinite height is again a module with no elements of infinite height.
\item[(ii)] A module with no elements of infinite height is pure in its $p$-adic completion.
\item[(iii)] If $T$ is a pure submodule of a complete module $M$ then the closure of $T$ is likewise pure.
\item[(iv)] A module with no elements of infinite height which is complete in its $p$-adic topology is the completion of a direct sum of cyclic modules.
\item[(v)] If $M$ is a $\Z_p$-module and $S$ is a pure submodule of $M$ with no elements of infinite height which is complete in its $p$-adic topology then $S$ is a direct summand of $M$.
\end{itemize}
\end{fact}

\begin{corollary}\label{corollary56}
If $A\in\K_p$ then $B\in\K_p$ and $B\ss A$ if and only if $B$ is a direct summand of $A$.
\end{corollary}

\begin{proof}
 If $A,B\in\K_p$ and $B\ss A$ then $B$ is a direct summand by (v) of Fact~\ref{fact55}. On the other hand if $B$ is a direct summand of $A$ then $B$ is a pure closed subgroup of $A$ and by (iv) of Fact~\ref{fact55} $B\in\K_p$. 
\end{proof}

Taking into account that a product of groups $A_i$ is complete in the $p$-adic topology if and only if every $A_i$ is, we may write our models in the form $\overline{\Z^{(\kappa)}}$. The backbone $\Z^{(\kappa)}$ of the model is what Fuchs \cite{fuchs} calls a $p$-basic subgroup:

\begin{definition}\label{definition57}
A $p$-basic subgroup $B$ of $A$ is a subgroup of $A$ satisfying the following three conditions:
\begin{itemize}
\item[(i)] $B$ is a direct sum of cyclic $p$-groups and infinite cyclic groups,
\item[(ii)] $B$ is $p$-pure in $A$ ($p^kB=B\cap p^kA$ for $k\in\N$),
\item[(iii)] $A/B$ is $p$-divisible ($p^kA/B=A/B$ for $k\in\N$).
\end{itemize}
If $B$ is a $p$-basic subgroup of $A$ then $B$ has a basis which is said to be a $p$-basis of $A$. This basis is $p$-independent, i.e., for every finite subsystem $a_1,\dots,a_m$ 
$$
n_1a_1+\cdots +n_ma_m\in pA \quad (n_ia_i\neq 0,n_i\in\Z)
$$
implies
$$
p\div n_i \quad (i=1,\dots,m).
$$
\end{definition}

\begin{fact}[\cite{fuchs}]\label{fact58}
\begin{itemize}
\item[(i)] A subgroup generated by a $p$-independent system in $A$ is $p$-pure in $A$.
\item[(ii)] Every $p$-independent system of $A$ can be expanded to a $p$-basis of $A$.
\item[(iii)] For a given prime $p$ all $p$-basic subgroups of a group are isomorphic.
\end{itemize}
\end{fact}

\begin{proposition}\label{proposition59}
For a given prime $p$, $\K_p$ is a MAEC with L\"owenheim-Skolem number $\aleph_0$.
\end{proposition}

\begin{proof}
 Both $\K_p$ and $\ss$ are closed under isomorphism. If $A\ss B$
then $A$ is a substructure of $B$ and $\ss$ is a partial order of
$\K$. For unions note that if $(A_i)$ is an increasing chain of groups in
$\K_p$ with $A_i$ pure in $A_j$ for $i\leq j$, then $\bigcup_iA_i$ is a
torsion-free $\Z_p$-module with no non-zero elements of infinite height.
By Fact~\ref{fact55} its completion is again a $\Z_p$-module with no elements of
infinite height and thus of the form $\overline{\Z_p^{(\kappa)}}$ and is
a model in $\K_p$.
Also, as each $A_i$ is pure in $\bigcup_iA_i$ which in turn by Fact~\ref{fact55} is pure in its completion, each $A_i$ is pure in $\overline{\bigcup_iA_i}$, and if for all $i$ $A_i\ss B\in\K_p$ then $\bigcup_iA_i$ is pure in $B$ and by Fact~\ref{fact55} so is $\overline{\bigcup_iA_i}$.

For the coherence axiom note that if $A\ss C$ then $A$ is a pure subgroup of any subgroup of $C$ that it is contained in. Thus if $B\ss C$ and $A\subset B$ we have $A\ss B$.

The L\"owenheim-Skolem number $LS^d(\K_p)$ is $\aleph_0$. To see this let $C$ be a subset of $A\in\K_p$. Clearly $\overline{\pcl{C}}$ is the smallest pure closed subgroup of $A$ containing $C$, so $C\subset\overline{\pcl{C}}\ss A$ and as $\pcl{C}$ has cardinality at most $\abs{C}+\aleph_0$ we are done. 
\end{proof}

In this example we only consider isometric isomorphisms so the $d^p$-metric reduces to the infimum-distance metric $d(p,q)=\inf\{d(a,b) : a\models p,b\models q\}$. Also almost summability trivially holds.

\begin{proposition}\label{proposition60}
The class $\K_p$ has the joint embedding and amalgamation properties.
\end{proposition}

\begin{proof}
Since direct sums of (disjoint) models are models this is clear by Corollary~\ref{corollary56}. 
\end{proof}

\begin{lemma}\label{lemma61}
If $A\ss B$ with a $p$-base of strictly smaller cardinality and if $f:A\to B$ is a $\K$-embedding, i.e., an embedding such that $f(A)\ss B$ then $f$ can be extended to an automorphism of $B$.
\end{lemma}

\begin{proof}
 Write $B=A\oplus B_1=f(A)\oplus B_2$ and note that by cardinality considerations $B_1$ and $B_2$ must be isomorphic. 
\end{proof}

By the above lemma any large enough model acts as a monster model and below we shall assume we work inside such a model $M\in\K_p$.

\begin{lemma}\label{lemma62}
If $A\subset M$ is a set, $a\in M$ an element and $a\notin\overline{\pcl{A}}$ then the Galois-type of $a$ over $A$, $t^g(a/A)$, is determined exactly by the distance of $a$ to $\pcl{A}$ and the set $A_a\subset\pcl{A}$ of closest elements, i.e. $A_a=\{b\in\pcl{A} : d(a,b)=d(a,\pcl{A})\}$.
\end{lemma}

\begin{proof}
 First note that $\pcl{A}$ (and thus also its closure) is fixed pointwise by any automorphism fixing $A$ pointwise. This is because in a torsion-free group the equation $nx=a$ has at most one solution. Thus also distances to elements within $\pcl{A}$ must be preserved.

Now if $a$ and $b$ have the same positive distance to $\pcl{A}$, say $p^{-k}$, and the same set of closest elements, choose one of these, say $c$ and write $a=p^ka'+c$, $b=p^kb'+c$, where $p\ndiv a',b'$. Let $I$ be a $p$-basis for $\pcl{A}$. Then $I\cup\{a'\}$ is $p$-independent: Let $a_1,\dots,a_m\in I$, $n_0,\dots,n_m\in\Z$ such that $p\div n_0a'+n_1a_1+\cdots+n_ma_m$. If for some $i\leq m$ $p\ndiv n_i$ we are left with three scenarios:
\begin{itemize}
\item[(i)] $p\div n_i$ for $1\leq i\leq m$ but $p\ndiv n_0$. But then $p\div a'$ contradicting the maximality of $k$ (in the distance of $a$ to $\pcl{A}$).
\item[(ii)] $p\div n_0$ but for some $1\leq i\leq m$ $p\ndiv n_i$. This contradicts $p$-independence of $I$.
\item[(iii)] $p\ndiv n_0$ and for some $1\leq i\leq m$ $p\ndiv n_i$. By removing terms we may assume $p\ndiv n_i$ for all $i\leq m$. As $p\ndiv n_0$ and $\pcl{A}$ is pure in the $\Z_p$-module $M$, there is $b\in\pcl{A}$ such that $n_0b=n_1a_1+\cdots n_ma_m$. Thus $p\div a'+b$, i.e. $p^{k+1}\div a-c+p^kb$ contradicting the maximality of $k$.
\end{itemize}
So $I\cup\{a'\}$ is $p$-independent and similarly $I\cup\{b'\}$. Thus we can construct a $\K$-embedding of $\pcl{Aa}$ into $M$ by fixing $I$ pointwise and mapping $a'$ to $b'$. By Lemma~\ref{lemma61} this extends to an automorphism of $M$. 
\end{proof}

Note that in the lemma the set of closest elements $A_a=\{b\in\pcl{A} : d(a,b)=d(a,\pcl{A})\}$ is a closed ball of radius $d(a,\pcl{A})$ so to check that two elements with the same distance to $\pcl{A}$  have the same type it is enough to show that their sets of closest elements intersect.

\begin{proposition}\label{proposition63}
The class $\K_p$ is homogeneous.
\end{proposition}

\begin{proof}
 Let $(a_i)_{i<\alpha}$ and $(b_i)_{i<\alpha}$ be sequences of elements in $M$ such that
$$t^g((a_{i_k})_{k<n}/\emptyset)=t^g((b_{i_k})_{k<n}/\emptyset)\quad{\rm for\ each}\ n<\omega.$$
Now define $f$ by $a_i\mapsto b_i$. As finite tuples of $(a_i)_{i<\alpha}$ and $(b_i)_{i<\alpha}$ have the same type, this induces a group isomorphism between $\pcl{(a_i)_{i<\alpha}}$ and $\pcl{(b_i)_{i<\alpha}}$ which naturally extends to their closures. Thus $f$ extends to a map $\overline{\pcl{(a_i)_{i<\alpha}}}$ to $\overline{\pcl{(b_i)_{i<\alpha}}}$ and as these are models the map further extends to an automorphism of $M$ by Lemma~\ref{lemma61}. 
\end{proof}

\begin{proposition}\label{proposition64}
The class $\K_p$ has the perturbation property, i.e., if $d^p(t^g(a/\emptyset),t^g(b/\emptyset))=0$ then $t^g(a/\emptyset)=t^g(b/\emptyset)$.
\end{proposition}

\begin{proof}
Let $(b_i)_{i<\omega}$ be a sequence of tuples in a large model $M$
such that $t^g(b_i/\emptyset)=t^g(b_j/\emptyset)$ for all $i,j<\omega$ and
assume $(b_i)_{i<\omega}$ converges to $b$. Denote $B_i=\pcl{b_i}$.
By assumption the mappings mapping $b_0$ to $b_i$ induce isomorphisms
$f_i:B_0\to B_i$. Now consider $B=\pcl{b}$. Define a map $f:B_0\to B$ by
letting $f(a)$, for $a\in B_0$, be the limit of $(f_i(a))_{i<\omega}$.
As $b_i\to b$ it is easy to see that linear combinations of $\cl{b_i}$
converge to the corresponding linear
combination of $\cl{b}$. Also each $f_i$ must preserve divisibility and if $(na_i')_{i<\omega}$ converges to some $a$ then $(a_i')_{i<\omega}$ must be convergent and its limit $a'$ must satisfy $na'=a$. Thus $f$ is a group isomorphism $B_0\to B$ and extends to the closures. These in turn are models, so $f$ extends to an automorphism of $M$ mapping $b_0\mapsto b$. 
\end{proof}

Note that as we only consider isometric mappings and $d^p$ thus coincides with the infimum-distance metric, completeness of type spaces ($d^p$-Cauchy sequences of types over $\emptyset$ have a limit) is just completeness of the model.

\begin{proposition}\label{proposition65}
The class $\K_p$ is $\omega$-$d^p$-stable, i.e., the set of types over a parameter set of cardinality (or density) $\aleph_0$ has density character $\aleph_0$ in the $d^p$-topology.
\end{proposition}

\begin{proof}
Let $A$ be a separable 
model of $\K_p$ and let $B=A\oplus\overline{\Z_p^{(\omega)}}$. It is
enough to show that all types over $A$ can be realised in $B$.
Since all types over $A$ can be realised in a separable strong extension
$C\succcurlyeq A$ we need to show that such a $C$ can be embedded over $A$
into $B$. Now as $A\ss C$, $A$ is a direct summand of $C$ so $C=A\oplus C'$
and $C'$ is either empty or of the form $\overline{\Z_p^{(\alpha)}}$.
Then $\alpha$ is
either finite or $\omega$ and thus $C'$ can be embedded into the complement of $A$ in $B$. Combining this with the identity map on $A$ we are done. 
\end{proof}

\begin{proposition}\label{proposition66}
For elements $a$, $a\da_AB$ if and only if $d(a,\pcl{A})=d(a,\pcl{AB})$.
\end{proposition}

\begin{proof}
 Assume $d(a,\pcl{AB})<d(a,\pcl{A})$ (which implies $d(a,\pcl{A})\neq 0$).
Choose an element $b\in\pcl{AB}$ such that $d(a,b)<d(a,\pcl{A})$
(and if $d(a,\pcl{AB})>0$ we can actually choose $b$ such that
$d(a,b)=d(a,\pcl{AB})$). Now let $b_A$ be a closest element to $b$ in
$\pcl{A}$. If $J_A$ is a $p$-basis for $\pcl{A}$ and $d(b,\pcl{A})=p^{-k}$
we can write $b=b_A+p^kb'_0$ where $J_A\cup\{b'_0\}$ is $p$-independent.
We can extend this to a $p$-basis $J_B$ for $\pcl{AB}$ and further to a
$p$-independent
sequence $J_B\cup\{b'_i\}_{0<i<\omega}$. Then define $b_i=b_A+p^kb'_i$. The $b_i$ form an $A$-indiscernible sequence, $b_0=b$ and for $i\neq j$ $d(b_i,b_j)=d(b_0,\pcl{A})=d(b,b_A)$. As $d(a,b)<d(a,b_A)$ we must have $d(a,b_A)=d(b,b_A)$. Now if $a'\models t^g(a/AB)$ then $d(a',b)=d(a,b)$ and we must have $d(a',b_1)=d(b_0,b_1)>d(a',b_0)$. So $t^g(a'/AB\cup{b'_i}_{0<i<\omega})$ splits strongly over $A$. This proves $a\nda_AB$.

For the other direction, assume $d(a,\pcl{A})=d(a,\pcl{AB})$. If this
distance is 0, then there is a countable $A'\subset A$ s.t.
$a\in\overline{\pcl{A'}}$. Then for any $B'\supseteq AB$, $t^g(a/B')$
does not split strongly over $A'$ as any $A'$-indiscernible sequence must
be $A'a$-indiscernible. If the distance is positive, let $A'\subset A$
be finite such that $d(a,\pcl{A'})=d(a,\pcl{A})$. Let $B'\supseteq AB$.
We need to find $b\models t^g(a/AB)$ such that $t^g(b/B')$ does not split
strongly over $A'$. As $d(a,\pcl{A'})=d(a,\pcl{AB})$
we can find $a_{A'}\in\pcl{A'}$ such that it is a closest element to $a$
in $\pcl{AB}$ and write $a=a_{A'}+a'$. Let $b=a_{A'}+b'$ where
$d(b',\pcl{B'})=\pnorm{a'}$. Then $b\models t^g(a/AB)$ and we prove
that $t^g(b/B')$ does not split strongly over $A'$: Let
$\{b_i:i<\omega\}\subset B'$ be $A'$-indiscernible. Then
$\overline{\pcl{A'b_0}}$ is a model and the map generated by fixing
$A'$ and mapping $b_0$ to $b_1$ is $\K$-elementary and by Lemma~\ref{lemma61} extends
to an automorphism of $\overline{\pcl{B'}}$. Now looking at $\overline{\pcl{B'}}$ inside any larger model containing $b'$ we see that we can extend the mapping to one fixing $b'$. Then we have a map showing that $t^g(bb_0/A')=t^g(bb_1/A')$. 
\end{proof}

\begin{corollary}\label{corollary67}
$\K_p$ is simple.
\end{corollary}

\begin{proof}
 Let $a$ be a finite tuple and $A$ a set. If $a$ is a single element then $a\da_AA$ by Proposition~\ref{proposition66}. If $a=a_1\dots a_n$, use induction on $n$ and Fact~\ref{fact2} (v). 
\end{proof}

\begin{lemma}\label{lemma68}
For single elements $a$, $R(a/A)=d(a,\pcl{A})$.
\end{lemma}

\begin{proof}
 If $a'\models t^g(a/A)$ then $d(a,a')\leq d(a,\pcl{A})$ so $a\da^{d(a,\pcl{A})}_AB$ for any $B$. On the other hand if $a\nda_AB$ and $a'\da_AB$ then $d(a,a')=d(a,\pcl{A})$ so $a\nda^{\e}_AB$ for all $\e<d(a,\pcl{A})$. 
\end{proof}

The last argument in the proof shows:

\begin{corollary}\label{corollary68}
For single elements $a$, if $a\da^{\e}_AB$ for some $\e<d(a,\pcl{A})$, then $a\da_AB$.
\end{corollary}

\begin{proposition}\label{proposition69}
In $\K_p$, for any set $A$ any type of a single element satisfies the assumptions of Corollary~\ref{corollary54}.
\end{proposition}

\begin{proof}
 Let $p=Lstp^w(a/A)$ where $a$ is a single element. Define $d=d(a,\pcl{A})$ and let $d^-$ be the largest distance smaller than $d$ (recall that the positive values of the metric are discrete). Choose any $\e$ within $d^-<\e<d$. Let $b\in D$ and $B\subseteq D$.

(i)$\Leftrightarrow$(iii) is clear by Lemma~\ref{lemma68} and Proposition~\ref{proposition66}.

(ii)$\Rightarrow$(i) is trivial.

(iv)$\Rightarrow$(ii): Assume (iv) holds for $c=b$, i.e., there exists $b'$ such that $b'\nda_Ab$ and $b\da^{>\e}_{AB}b'$. Now if $b\da^{>\e}_AB$, by Corollary~\ref{corollary68} and the choice of $\e$ $b\da_AB$. Then by Corollary~\ref{corollary43} $b\da^{>\e}_ABb'$ and again by Corollary~\ref{corollary68} $b\da_Ab'$, a contradiction.

(i)$\Rightarrow$(iv): Assume $b\nda_AB$. Then $d(b,\pcl{AB})\leq d^-$ so there is $b'\in\pcl{AB}$ satisfying $d(b,b')\leq d^-$. Then $b'\nda_Ab$ and $\pcl{ABb'}=\pcl{AB}$ so any element $c\in D$ must satisfy $c\da_{AB}b'$ and thus $c\da^{>e}_{AB}b'$. 
\end{proof}

Thus the $\nda_A$-equivalence classes of a type (over $A$) form a pregeometry.
Note that we really need to look at equivalence classes to get a pregeometry.
If we define the closure simply by $cl(B)=\{a:a\nda_AB\}$ the property
$cl(cl(B))=cl(B)$ fails. This can be seen by considering $p$-independent
elements $b_i$ and letting, e.g., $A=\emptyset$, $p$ be the type of any
element of length $1$, $B=\{b_1-pb_2\}$, $c=pb_0+b_1$ and $a=b_0+b_2$.
Then $c\in cl(B)$, $a\in cl(Bc)$ but $a\notin cl(B)$. This reflects the way the structure theorem for (nice)
Abelian groups looks at the Ulm invariants, i.e., the dimensions of $p^\a G/p^{\a+1}G$, considered as  vector spaces over the integers mod $p$. 
Note that when $A=\emptyset$ and $a,b\in G-pG$, they are in the same $\nda$-equivalence class if and only if, when $G/pG$ is considered as a vector space over $\Z/p\Z$, the cosets of $a$ and $b$ span the same linear subspace, $(\Z/p\Z)a/pG=(\Z/p\Z)b/pG$.

\end{document}